\documentclass[12pt,twoside,reqno,leqno]{amsart}
\usepackage{url}
\usepackage{amsmath}
\usepackage{amssymb}
\usepackage{amsfonts}
\usepackage{amsthm}
\usepackage[mathscr]{eucal}
\usepackage{amscd}
\usepackage{amsmath}
\usepackage{latexsym}
\usepackage[dvips]{graphics}
\usepackage{graphics,epstopdf}
\usepackage[pdftex]{graphicx}
\theoremstyle{plain}

\newtheorem{prop}{Proposition}
\theoremstyle{definition}

 \theoremstyle{remark}
 \newtheorem{rem}{Remark}
\newtheorem{example}{Example}

\allowdisplaybreaks
\def\DJ{\leavevmode\setbox0=\hbox{D}\kern0pt\rlap
{\kern.04em\raise.188\ht0\hbox{-}}D}
\begin{document}
\title[Unified Multi-Tupled Fixed Point Theorems]
{Unified Multi-Tupled Fixed Point Theorems involving Monotone
Property in Ordered Metric Spaces}
\author[Alam, Imdad and Radenovi$\acute{\rm c}$]{Aftab Alam$^1$, Mohammad Imdad$^2$ and Stojan Radenovi$\acute{\rm c}$$^3$}
\maketitle
\begin{center}
{\footnotesize $^{1,2}$Department of Mathematics, Aligarh Muslim University, Aligarh-202002, Uttar Pradesh, India.\\
$^3$University of Belgrade, Faculty of Mechanical Engineering, Kraljice Marije 16, 11120 Beograd, Serbia.\\
Email addresses: aafu.amu@gmail.com,  mhimdad@yahoo.co.in, radens@beotel.rs.\\}
\end{center}
{\footnotesize{\noindent {\bf Abstract.}} In this paper, we
introduce a generalized notion of monotone property and prove some
results regarding existence and uniqueness of multi-tupled fixed
points for nonlinear contraction mappings satisfying monotone
property in ordered complete metric spaces. Our results unify
several classical and well known $n$-tupled (including coupled,
tripled and quadruple ones) fixed point results existing in
literature. \vskip0.5cm \noindent {\bf Keywords}:  {\it ICU}
property; Monotone property, $\ast$-fixed points,
$\varphi$-contractions.

}
\section{Introduction}
Throughout this manuscript, the following symbols and notations are
involved:
\begin{enumerate}
\item [{(1)}] As usual,$(X,d)$, $(X,\preceq)$ and $(X,d,\preceq)$ are termed as metric space, ordered set and ordered metric space, wherein $X$ stands for a
nonempty set $X$, $d$ for a metric on $X$ and $\preceq$ for a
partial order on $X$. Moreover, if the metric space $(X,d)$ is
complete, then $(X,d,\preceq)$ is termed as ordered complete metric
space.
\item [{(2)}] $\succeq$ denotes dual partial order of
$\preceq$ ($i.e$ $x\succeq y$ means $y\preceq x$).
\item [{(3)}] $\mathbb{N}$ and $\mathbb{N}_0$ stands for the sets of positive and nonnegative
integers respectively ($i.e.$ $\mathbb{N}_0=\mathbb{N}\cup \{0\}$).
\item [{(4)}] $n$ stands for a fixed natural number greater than
 1, while $m,l \in \mathbb{N}_0$.
\item [{(5)}] $I_n$ denotes the set $\{1,2,...,n\}$ and we use $i,j,k\in I_n.$
\item [{(6)}] For a nonempty set $X$, $X^n$ denotes the cartesian
product of $n$ identical copies of $X$, $i.e.$, $X^n:=X\times X
\times\stackrel{(n)}{...}\times X$. We call $X^n$ the
$n$-dimensional product set induced by $X$.
\item [{(7)}] A sequence in $X$ is denoted by $\{x^{(m)}\}$ and a
sequence in $X^n$ is denoted by $\{{\rm U}^{(m)}\}$ where
U$^{(m)}=(x^{(m)}_1,x^{(m)}_2,...,x^{(m)}_n)$ such that for each
$i\in I_n$, $\{x^{(m)}_i\}$ is a sequence in $X$.
\end{enumerate}

Starting from the Bhaskar-Lakshmikantham coupled fixed point theorem
($cf.$ \cite{C1}), the branch of multi-tupled fixed point theory in
ordered metric spaces is progressed in high speed during only one
decade. Then, coupled fixed point theorems are extended upto higher
dimensional product set by appearing tripled, quadrupled and
$n$-tupled fixed point theorems. Here it can be highlighted that
extension of coupled fixed point upto higher dimensional product set
is not unique. It is defined by various author in different way.
From recent year, some authors paid attention to unify the different
types of multi-tupled fixed points. A first attempt of this kind was
given by Berzig and Samet \cite{HD1}, wherein the authors defined a
unified notion of $n$-tupled fixed point by using 2$n$ mappings from
$I_n$ to $I_n$. Later, Rold$\acute{\rm a}$n $et\;al.$ \cite{MD1}
extended the notion of $n$-tupled fixed point of Berzig and Samet
\cite{HD1} by introducing the notion of $\Upsilon$-fixed point based
on $n$ mappings from $I_n$ to $I_n$. Most, recently, Alam $et\;al.$
\cite{unif1} modified the notion of $\Upsilon$-fixed point by
introducing the notion of $\ast$-fixed point depending only one
mapping After appearance of multi-tupled fixed points, some authors
paid attention to unify the different types of multi-tupled fixed
points. A first attempt of this kind was given by Berzig and Samet
\cite{HD1}, wherein the authors defined a unified notion of
$n$-tupled fixed point by using 2$n$ mappings from $I_n$ to $I_n$.
Later, Rold$\acute{\rm a}$n $et\;al.$ \cite{MD1} extended the notion
of $n$-tupled fixed point of Berzig and Samet \cite{HD1} by
introduced the notion of $\Upsilon$-fixed point based on $n$
mappings from $I_n$ to $I_n$. Most, recently, Alam $et\;al.$
\cite{unif1} modified the notion of $\Upsilon$-fixed point by
introducing the notion of $\ast$-fixed point depending on a binary
operation $\ast$ on $I_n$.\\

One of the common properties of multi-tupled fixed point theory in
the context of ordered metric spaces is that the mapping
$F:X^n\rightarrow X$ satisfies mixed monotone property (for
instance, see \cite{C2}-\cite{Q2}). In order to avoid the mixed
monotone property in such results, authors in \cite{MT1}-\cite{MN}
utilized the notion of monotone property.\\

The aim of this paper, is to extend the notion of monotone property
for the mapping $F:X^n\rightarrow X$ and utilizing this, to prove
some existence and uniqueness results on $\ast$-coincidence points
under $\varphi$-contractions due to Boyd and Wong \cite{B2}.\\

\section{Monotone property in ordered metric spaces}

In this section, we recall some initial results regarding monotone
property and then by motivating this, we introduced
generalized notion of monotone property. But before starting these results, we summarize some preliminaries used in such results. \\

Given two self-mappings $f,g$ defined on an ordered set
$(X,\preceq)$, we say that $f$ is $g$-increasing if for any $x,y\in
X$, $g(x)\preceq g(y)$ implies $f(x)\preceq f(y)$ ($cf.$
\cite{PGF3}). As per standard practice, we can defined the notions
of increasing, decreasing, monotone, bounded above and bounded below
sequences besides bounds (upper as well as lower) of a sequence in
an ordered set $(X,\preceq)$, which on the set of real numbers with
natural ordering coincide with their usual senses.\\

Let $(X,d,\preceq)$ be an ordered metric space and $\{x_n\}$ a
sequence in $X$. We adopt the following notations.
\begin{enumerate}
\item[{(i)}] If $\{x_n\}$ is increasing and $x_n\stackrel{d}{\longrightarrow} x$, then we denote it symbolically by
$x_n\uparrow x.$
\item[{(ii)}] If $\{x_n\}$ is decreasing and $x_n\stackrel{d}{\longrightarrow} x$, then we denote it symbolically
by $x_n\downarrow x.$
\item[{(iii)}] If $\{x_n\}$ is monotone and $x_n\stackrel{d}{\longrightarrow} x$, then we denote it symbolically
by $x_n\uparrow\downarrow x.$\\
\end{enumerate}

Alam $et\;al.$ \cite{PGF13} formulated the following notions by
using certain properties on ordered metric space (in order to avoid
the necessity of continuity requirement on underlying mapping)
utilized by earlier authors especially from \cite{PF2,PGF3,C1,C2} besides some other ones.\\

\noindent\textbf{Definition 1} \cite{PGF13}. Let $(X,d,\preceq)$ be
an ordered metric space and $g$ a self-mapping on $X$. We say that
\begin{enumerate}
\item[{(i)}] $(X,d,\preceq)$ has {\it
g-ICU}\;(increasing-convergence-upper bound) property if $g$-image
of every increasing convergent sequence $\{x_n\}$ in $X$ is bounded
above by $g$-image of its limit (as an upper bound), $i.e.,$
$$x_n\uparrow x \Rightarrow g(x_n)\preceq g(x)\;\;\forall~ n\in \mathbb{N}_0,$$
\item[{(ii)}] $(X,d,\preceq)$ has {\it
g-DCL}\;(decreasing-convergence-lower bound) property if $g$-image
of every decreasing convergent sequence $\{x_n\}$ in $X$ is bounded
below by $g$-image of its limit (as a lower bound), $i.e.,$
$$x_n\downarrow x \Rightarrow g(x_n)\succeq g(x)\;\;\forall~ n\in \mathbb{N}_0\;{\rm and}$$
\item[{(iii)}] $(X,d,\preceq)$ has {\it g-MCB}\;(monotone-convergence-boundedness)
property if $X$ has both {\it g-ICU} as well as {\it g-DCL}
property.
\end{enumerate}

Notice that under the restriction $g=I,$ the identity mapping on
$X,$ the notions of {\it g-ICU} property, {\it g-DCL} property and {\it g-MCB} property are respectively called {\it ICU} property, {\it DCL} property and {\it MCB} property.\\

\noindent\textbf{Definition 2} (Borcut \cite{MT1}). Let
$(X,\preceq)$ be an ordered set and $F: X^{3} \rightarrow X$ a
mappings. We say that $F$ is monotone if $F(x,y,z)$ is increasing in
$x,y,z$, $i.e.$, for any $x,y,z\in X,$
$$x_{1}, x_{2}\in X, x_{1}\preceq x_{2}\Rightarrow F(x_{1},y,z)\preceq F(x_{2},y,z),$$
$$y_{1}, y_{2}\in X, y_{1}\preceq y_{2}\Rightarrow F(x,y_{1},z)\preceq F(x,y_{2},z),$$
$$z_{1}, z_{2}\in X, z_{1}\preceq z_{2}\Rightarrow F(x,y,z_1)\preceq F(x,y,z_2).$$

\noindent\textbf{Definition 3} (Borcut \cite{MT1}). Let $X$ be a
nonempty set and $F:X^{3} \rightarrow X$ a mapping. An element
$(x,y,z)\in X^3$ is called a tripled fixed point of $F$ if
$$F((x,y,z)=x,\;F(y,x,z)=y,\;F(z,y,x)=z.$$

Notice that this concept of tripled fixed point is essentially
different to that
of Berinde and Borcut \cite{T1}.\\

Utilizing the above notions, Borcut \cite{MT1} proved the following
variant of Berinde-Borcut tripled fixed point theorem \cite{T1}
without mixed
monotone property.\\

 \noindent{\bf Theorem 1.} (Borcut \cite{MT1}). Let
$(X,d,\preceq)$ be an ordered complete metric space and
$F:X^{3}\rightarrow X$ a mapping. Suppose that the following
conditions hold:
\begin{enumerate}
\item [{(i)}] $F$ is monotone,
\item [{(ii)}] either $F$ is continuous or $(X,d,\preceq)$ has {\it ICU} property,
\item [{(iii)}] there exist $x^{(0)},y^{(0)},z^{(0)}
\in X$ such that
$$ x^{(0)} \preceq
F(x^{(0)},y^{(0)},z^{(0)}), y^{(0)}\preceq
F(y^{(0)},x^{(0)},z^{(0)})\;{\rm and }\; z^{(0)} \preceq
F(z^{(0)},y^{(0)},x^{(0)}),$$
\item [{(iv)}] there exist $\alpha,\beta,\gamma\in [0,1)$ with $\alpha+\beta+\gamma<1$ such that
$$d(F(x,y,z),F(u,v,w))\leq\alpha d(x,u)+\beta d(y,v)+ \gamma d(z,w)$$
for all $x,y,z,u,v,w\in X$ with $x\preceq u$, $y\preceq v$ and
$z\preceq w$.
\end{enumerate}
Then $F$ has a tripled fixed point (in the sense of Borcut
\cite{MT1}), $i.e.$, there exist $x,y,z\in X$ such that
$F(x,y,z)=x$, $F(y,x,z)=y$ and $F(z,y,x)=z$.\\

Later, Borcut \cite{MT2} generalized above concept for a pair of
mappings $F:X^{3}\rightarrow X$ and $g:X\to X$ as follows:\\

\noindent\textbf{Definition 4} (Borcut \cite{MT2}). Let
$(X,\preceq)$ be an ordered set and $F:X^{3}\rightarrow X$ and
$g:X\to X $ two mappings. We say that $F$ has $g$-monotone property
if $F(x,y,z)$ is $g$-increasing in $x,y,z$, $i.e.$, for any
$x,y,z\in X,$
$$x_{1}, x_{2}\in X, g(x_{1})\preceq g(x_{2})\Rightarrow F(x_{1},y,z)\preceq F(x_{2},y,z),$$
$$y_{1}, y_{2}\in X, g(y_{1})\preceq g(y_{2})\Rightarrow F(x,y_{1},z)\preceq F(x,y_{2},z),$$
$$z_{1}, z_{2}\in X, g(z_{1})\preceq g(z_{2})\Rightarrow F(x,y,z_1)\preceq F(x,y,z_2).$$

\noindent\textbf{Definition 5} (Borcut \cite{MT2}). Let $X$ be a
nonempty set and $F:X^{3}\rightarrow X$ and $g:X\to X $ two
mappings. An element $(x,y,z)\in X^3$ is called a tripled
coincidence point of $F$ and $g$ if
$$F((x,y,z)=g(x),\;F(y,x,z)=g(y),\;F(z,y,x)=g(z).$$

\noindent\textbf{Definition 6} (Borcut \cite{MT2}). Let $X$ be a
nonempty set and $F: X^{3} \rightarrow X$ and $g: X\rightarrow X$
two mappings. We say that the pair $(F,g)$ is commuting if
$$g(F(x,y,z))=F(gx,gy,gz)\;\;\forall ~x,y,z\in
X.$$

\noindent\textbf{Definition 7} (Lakshmikantham and \'{C}iri\'{c} \cite{C2}). We denote by $\Phi$ the family of functions $\varphi : [0,\infty)\to [0,\infty)$ satisfying\\
\indent\hspace{1cm} $(a)$ $\varphi(t)<t$ for each $t>0$ ,\\
\indent\hspace{1cm} $(b)$ $\lim\limits_{r\to t^+}\varphi(r)<t$ for
each $t>0$.\\

Borcut \cite{MT2} proved the following tripled coincidence theorem
for nonlinear contraction satisfying $g$-monotone property.\\

\noindent{\bf Theorem 2.} (Borcut \cite{MT2}).  Let $(X,d,\preceq)$
be an ordered complete metric space and $F:X^{3}\rightarrow X$ and
$g:X\to X $ two mappings. Suppose that the following conditions
hold:
\begin{enumerate}
\item [{(i)}] $F(X^3)\subseteq g(X)$,
\item [{(ii)}] $F$ has $g$-monotone property,
\item [{(iii)}] $(F,g)$ is commuting pair,
\item [{(iv)}] $g$ is continuous,
\item [{(v)}] either $F$ is continuous or $(X,d,\preceq)$ has $g$-{\it ICU}
property,
\item [{(vi)}] there exist $x^{(0)},y^{(0)},z^{(0)}
\in X$ such that
$$ g(x^{(0)}) \preceq
F(x^{(0)},y^{(0)},z^{(0)}), g(y^{(0)})\preceq
F(y^{(0)},x^{(0)},z^{(0)})\;{\rm and }\; g(z^{(0)}) \preceq
F(z^{(0)},y^{(0)},x^{(0)}),$$
\item [{(vii)}] there exists $\varphi\in \Phi$ such that
$$d(F(x,y,z),F(u,v,w))\leq\varphi(\max\{d(gx,gu),d(gy,gv),d(gz,gw)\})$$
for all $x,y,z,u,v,w\in X$ with $g(x)\preceq g(u)$, $g(y)\preceq
g(v)$ and $g(z)\preceq g(w)$.
\end{enumerate}

Then $F$ and $g$ have a tripled coincidence point(in the sense of
Borcut \cite{MT2}), $i.e.$, there exist $x,y,z\in X$ such that
$F(x,y,z)=g(x)$, $F(y,x,z)=g(y)$ and $F(z,y,x)=g(z)$.\\

Recently, Radenovi$\acute{\rm c}$ \cite{MC1,MC2} discussed
bi-dimensional variant of monotonicity.\\

\noindent\textbf{Definition 8} (Radenovi$\acute{\rm c}$ \cite{MC1}).
Let $(X,\preceq)$ be an ordered set and $F: X^{2} \rightarrow X$ a
mappings. We say that $F$ is monotone if $F(x,y)$ is increasing in
$x,y$, $i.e.$, for any $x,y\in X,$
$$x_{1}, x_{2}\in X, x_{1}\preceq x_{2}\Rightarrow F(x_{1},y)\preceq F(x_{2},y),$$
$$y_{1}, y_{2}\in X, y_{1}\preceq y_{2}\Rightarrow F(x,y_{1})\preceq F(x,y_{2}).$$

\noindent\textbf{Definition 9} (Guo and Lakshmikantham \cite{C04}).
Let $X$ be a nonempty set and $F:X^{2} \rightarrow X$ a mapping. An
element $(x,y)\in X$ is called a coupled fixed point of $F$ if
$$F(x,y)=x\;{\rm and}\;F(y,x)=y.$$

The following result is the formulation of Bhaskar-Lakshmikantham
coupled fixed point theorem (\cite{C1}) for monotone mappings.\\

\noindent{\bf Theorem 3.} (Radenovi$\acute{\rm c}$ \cite{MC1}). Let
$(X,d,\preceq)$ be an ordered complete metric space and
$F:X^{2}\rightarrow X$ a mapping. Suppose that the following
conditions hold:

\begin{enumerate}
\item [{(i)}] $F$ is monotone,
\item [{(ii)}] either $F$ is continuous or $(X,d,\preceq)$ has {\it MCB} property,
\item [{(iii)}] there exist $x^{(0)},y^{(0)}\in X$ such that $x^{(0)}\preceq
F(x^{(0)},y^{(0)})\;{\rm and }\; y^{(0)} \preceq
F(y^{(0)},x^{(0)}),$
\item [{(iv)}] there exists $\alpha\in [0,1)$ such that
$$d(F(x,y),F(u,v))\leq\frac{\alpha}{2} [d(x,u)+ d(y,v)]$$
for all $x,y,u,v\in X$ with $x\preceq u$ and $y\preceq v$.
\end{enumerate}
Then $F$ has a coupled fixed point.\\

Here it can be pointed out that merely {\it ICU} property can serve our purpose instead of  {\it MCB} property.\\

\noindent\textbf{Definition 10} (Radenovi$\acute{\rm c}$
\cite{MC2}). Let $(X,\preceq)$ be an ordered set and $F: X^{2}
\rightarrow X$ and $g:X\rightarrow X$ two mappings. We say that $F$
has $g$-monotone property if $F(x,y)$ is $g$-increasing in $x,y$,
$i.e.$, for any $x,y\in X,$
$$x_{1}, x_{2}\in X, g(x_{1})\preceq g(x_{2})\Rightarrow F(x_{1},y)\preceq F(x_{2},y),$$
$$y_{1}, y_{2}\in X, g(y_{1})\preceq g(y_{2})\Rightarrow F(x,y_{1})\preceq F(x,y_{2}).$$

\noindent\textbf{Definition 11} (Lakshmikantham and \'{C}iri\'{c}
\cite{C2}). Let $X$ be a nonempty set and $F: X^{2} \rightarrow X$
and $g: X\rightarrow X$ two mappings. An element $(x,y)\in X^2$ is
called a coupled coincidence point of mappings $F$ and $g$ if
$$F(x,y)=g(x),\; F(y,x)=g(y).$$

\noindent\textbf{Definition 12} (Choudhury and Kundu \cite{C3}). Let
$(X,d)$ be a metric space and $F: X^{2} \rightarrow X$ and $g:
X\rightarrow X$ two mappings. We say that the pair $(F,g)$ is
compatible if
$$\lim\limits_{n\to \infty}d(gF(x_n,y_n),F(gx_n,gy_n))=0$$
and\indent\hspace{3.5cm} $\lim\limits_{n\to \infty}d(gF(y_n,x_n),F(gy_n,gx_n))=0$,\\
whenever $\{x_n\}$ and $\{y_n\}$ are sequences in $X$ such that
$$\lim\limits_{n\to \infty}F(x_n,y_n)=\lim\limits_{n\to \infty}g(x_n)\;{\rm and}\; \lim\limits_{n\to \infty}F(y_n,x_n)=\lim\limits_{n\to \infty}g(y_n).$$

\noindent\textbf{Theorem 4} (Radenovi$\acute{\rm c}$ \cite{MC2}).
Let $(X,d,\preceq)$ be an ordered metric space and $F:X^2\rightarrow
X$ and $g:X\rightarrow X$ two mappings. Assume that there exists
$\varphi\in \Phi$ such that
$$\max\{d(F(x,y),F(u,v)),d(F(y,x),F(v,u))\}\leq\varphi(\max\{d(gx,gu),d(gy,gv)\})$$
for all $x,y,u,v\in X$ with $g(x)\preceq g(u)$ and $g(y)\preceq
g(v)$ or $g(x)\succeq g(u)$ and $g(y)\succeq
g(v)$.\\
If the following conditions hold:\\
\indent\hspace{0.5mm}(i) $F(X^2)\subseteq g(X)$,\\
\indent\hspace{0.5mm}(ii) $F$ has the $g$-monotone property,\\
\indent\hspace{0.5mm}(iii) there exist $x^{(0)},y^{(0)}\in X$ such
that $$g(x^{(0)})\preceq F(x^{(0)},y^{(0)})\;{\rm and }\; g(y^{(0)})
\preceq F(y^{(0)},x^{(0)})$$ or $$g(x^{(0)})\succeq
F(x^{(0)},y^{(0)})\;{\rm and }\; g(y^{(0)}) \succeq
F(y^{(0)},x^{(0)})$$
\indent\hspace{0.5mm}(iv) $F$ and $g$ are continuous and compatible and $(X,d)$ is complete, or\\
\indent\hspace{0.5mm}(v) $(X,d,\preceq)$ has {\it MCB} property and one of $F(X^2)$ or $g(X)$ is complete.\\
Then $F$ and $g$ have a coupled coincidence point.\\

Now, we define generalized notions of monotone property as follows:\\

\noindent\textbf{Definition 13}. Let $(X,\preceq)$ be an ordered set
and $F: X^{n} \rightarrow X$ and $g: X\rightarrow X$ two mappings.
We say that $F$ has argumentwise $g$-monotone property if $F$ is
$g$-increasing in each of its arguments, $i.e.$, for any $x_1,x_2,...,x_n\in X$ and for each $i\in I_n$, \\
\indent\hspace{1cm}$\underline{x}_i,\overline{x}_i\in X$,
$g(\underline{x}_i)\preceq g(\overline{x}_i)$
$$\Rightarrow F(x_{1},x_{2},...,x_{i-1},\underline{x}_i,x_{i-1},...,x_{n})\preceq
F(x_{1},x_{2},...,x_{i-1},\overline{x}_i,x_{i-1},...,x_{n})$$

\noindent\textbf{Definition 14}. Let $(X,\preceq)$ be an ordered set
and $F: X^{n} \rightarrow X$ and $g: X\rightarrow X$ two mappings.
We say that $F$ has $g$-monotone property if for any
$x_1,x_2,...,x_n,y_1,y_2,...,y_n\in X,$
$$g(x_{1})\preceq g(y_1),g(x_{2})\preceq g(y_2),...,g(x_{n})\preceq g(y_n)$$
$$ \Rightarrow F(x_1,x_2,...,x_n)\preceq F(y_1,y_2,...,y_n).$$

On particularizing with $g=I$, the identity mapping on $X$, the
notions employed in Definitions 13 and 14 are, respectively, called
argumentwise monotone property
and monotone property.\\

Notice that the notions of monotone mappings followed in Definitions
2 and 8 and monotone property followed in Definitions 4 and 10 are
same as the notion of argumentwise monotone property followed in
Definition 13 but different from monotone property followed in
Definition 14. Henceforth, coherently with Definition 13, we prefer
calling these notions employing the term argumentwise monotone
property.

It is clear if $F$ has argumentwise monotone property (resp.
argumentwise $g$-monotone property) then it also has monotone
property (resp. $g$-monotone property).\\

\section{Coincidence Theorems in Ordered Metric Spaces}

In order to obtain multi-tupled fixed point theorems from
corresponding coincidence theorems, we indicate two recently proved
coincidence theorems. We know that in the context of order-theoretic
metrical fixed point theory, contractivity condition is compatible
with underlying partial order. Recently, Alam $et\;al.$ \cite{PGF17}
and Alam and  Imdad \cite{PGF18} defined compatibility of another
metrical notions with partial order by introducing their O,
$\overline{\rm
O}$ and $\overline{\rm O}$ analogous, which are summarized as follows:\\

\noindent\textbf{Definition 15} (Alam $et\;al.$ \cite{PGF18}). Let
$(X,d,\preceq)$ be called an ordered metric space. A nonempty subset
$Y$ of $X$ is
called a subspace of $X$ if $Y$ itself is an ordered metric space equipped with the metric $d_{Y}$ and partial order $\preceq_{Y}$ defined by:\\
$$d_{Y}(x,y)=d(x,y)\;\forall~x,y\in Y$$ and
$$x \preceq_{Y} y\Leftrightarrow x\preceq y\;\forall~x,y\in Y.$$

\noindent\textbf{Definition 16} (Alam $et\;al.$ \cite{PGF17}). An
ordered metric space $(X,d,\preceq)$ is called
\begin{enumerate}
\item[{(i)}] $\overline{\rm O}$-complete if every increasing Cauchy sequence in $X$ converges,
\item[{(ii)}] $\underline{\rm O}$-complete if every decreasing Cauchy sequence in $X$
converges and
\item[{(iii)}] O-complete if every monotone Cauchy sequence in $X$ converges.
\end{enumerate}
\begin{rem} \cite{PGF17} In an ordered
metric space, completeness $\Rightarrow$ O-completeness
$\Rightarrow$ $\overline{\rm O}$-completeness as well as
$\underline{\rm O}$-completeness.\end{rem}

\noindent\textbf{Definition 17} (Alam and  Imdad \cite{PGF18}). Let
$(X,d,\preceq)$ be an ordered metric space. A subset $E$ of $X$ is
called
\begin{enumerate}
\item[{(i)}] $\overline{\rm O}$-closed if for any sequence $\{x_n\}\subset E$,
$$x_n\uparrow x\Rightarrow x\in E,$$
\item[{(ii)}] $\underline{\rm O}$-closed if for any sequence $\{x_n\}\subset E$,
$$x_n\downarrow x\Rightarrow x\in E\;{\rm and}$$
\item[{(iii)}] O-closed if for any sequence $\{x_n\}\subset E$,
$$x_n \uparrow\downarrow x\Rightarrow x\in E.$$
\end{enumerate}

\begin{rem} \cite{PGF18} In an ordered metric space, closedness $\Rightarrow$ O-closedness
$\Rightarrow$ $\overline{\rm O}$-closedness as well as
$\underline{\rm O}$-closedness.\end{rem}

\begin{prop} \cite{PGF18}. Let $(X,d,\preceq)$
be an ordered metric space and $Y$ be a subspace of $X$.\\
\indent\hspace{5mm}(i) If $X$ is $\overline{\rm O}$-complete then
$Y$ is $\overline{\rm O}$-closed iff $Y$ is $\overline{\rm
O}$-complete.\\
\indent\hspace{5mm}(ii) If $X$ is $\underline{\rm O}$-complete then
$Y$ is $\underline{\rm O}$-closed iff $Y$ is $\overline{\rm
O}$-complete.\\
\indent\hspace{5mm}(iii) If $X$ is {\rm O}-complete then $Y$ is {\rm
O}-closed iff $Y$ is {\rm O}-complete.\end{prop}

\noindent\textbf{Definition 18} (Alam $et\;al.$ \cite{PGF17}). Let
$(X,d,\preceq)$ be an ordered metric space, $f:X\rightarrow X$ a
mapping and $x\in X$. Then $f$ is called
\begin{enumerate}
\item[{(i)}] $\overline{\rm O}$-continuous at $x$ if for any sequence $\{x_n\}\subset X$,
$$x_n\uparrow x\Rightarrow f(x_n)\stackrel{d}{\longrightarrow} f(x),$$
\item[{(ii)}] $\underline{\rm O}$-continuous at $x$ if for any sequence $\{x_n\}\subset X$,
$$x_n\downarrow x\Rightarrow f(x_n)\stackrel{d}{\longrightarrow} f(x)\;{\rm and}$$
\item[{(iii)}] O-continuous at $x$ if for any sequence $\{x_n\}\subset X$,
$$x_n \uparrow\downarrow x\Rightarrow f(x_n)\stackrel{d}{\longrightarrow} f(x).$$
\end{enumerate}
Moreover, $f$ is called O-continuous (resp. $\overline{\rm
O}$-continuous, $\underline{\rm O}$-continuous) if it is
O-continuous (resp. $\overline{\rm O}$-continuous, $\underline{\rm
O}$-continuous) at each point of $X$.
\begin{rem} \cite{PGF17} In an ordered
metric space, continuity $\Rightarrow$ O-continuity $\Rightarrow$
$\overline{\rm O}$-continuity as well as $\underline{\rm
O}$-continuity .\end{rem}

\noindent\textbf{Definition 19} (Alam $et\;al.$ \cite{PGF17}). Let
$(X,d,\preceq)$ be an ordered metric space, $f$ and $g$ two
self-mappings on $X$ and $x\in X$. Then $f$ is called
\begin{enumerate}
\item[{(i)}] $(g,\overline{\rm O})$-continuous at $x$ if for any sequence $\{x_n\}\subset X$,
$$g(x_n)\uparrow g(x)\Rightarrow f(x_n)\stackrel{d}{\longrightarrow} f(x),$$
\item[{(ii)}] $(g,\underline{\rm O})$-continuous at $x$ if for any sequence $\{x_n\}\subset X$,
$$g(x_n)\downarrow g(x)\Rightarrow f(x_n)\stackrel{d}{\longrightarrow} f(x)\;{\rm and}$$
\item[{(iii)}] $(g,{\rm O})$-continuous at $x$ if for any sequence $\{x_n\}\subset X$,
$$g(x_n) \uparrow\downarrow g(x)\Rightarrow f(x_n)\stackrel{d}{\longrightarrow} f(x).$$
\end{enumerate}
Moreover, $f$ is called $(g,{\rm O})$-continuous (resp.
$(g,\overline{\rm O})$-continuous, $(g,\underline{\rm
O})$-continuous) if it is $(g,{\rm O})$-continuous (resp.
$(g,\overline{\rm O})$-continuous, $(g,\underline{\rm
O})$-continuous) at each point of $X$.\\

Notice that on setting $g=I$ (the identity mapping on $X$),
Definition 16 reduces to Definition 15.\\

\noindent\textbf{Definition 20} (Alam $et\;al.$ \cite{PGF17}). Let
$(X,d,\preceq)$ be an ordered metric space and $f$ and $g$ two
self-mappings on $X$. We say that $f$ and $g$ are
\begin{enumerate}
\item[{(i)}] $\overline{\rm O}$-compatible if for any sequence $\{x_n\}\subset
X$ and for any $z\in X$,
$$g(x_n)\uparrow z\;{\rm and}\;f(x_n)\uparrow z\Rightarrow\lim\limits_{n\to \infty}d(gfx_n,fgx_n)=0,$$
\item[{(ii)}] $\underline{\rm O}$-compatible if for any sequence $\{x_n\}\subset
X$ and for any $z\in X$,
$$g(x_n)\downarrow z\;{\rm and}\;f(x_n)\downarrow z\Rightarrow\lim\limits_{n\to \infty}d(gfx_n,fgx_n)=0\;{\rm and}$$
\item[{(iii)}] O-compatible if for any sequence $\{x_n\}\subset
X$ and for any $z\in X$,
$$g(x_n)\uparrow\downarrow z\;{\rm and}\;f(x_n)\uparrow\downarrow z\Rightarrow\lim\limits_{n\to \infty}d(gfx_n,fgx_n)=0.$$
\end{enumerate}

The following family of control functions is indicated in  Boyd and
Wong \cite{B2} but was later used in
Jotic \cite{B4}.\\

$$\Omega=\Big\{\varphi:[0,\infty)\to [0,\infty):\varphi(t)<t\;{\rm for~each}\; t>0\;{\rm and~}\limsup\limits_{r\to t^+}\varphi(r)<t \;{\rm for~each~t>0}\Big\}.$$
Recently, Alam $et\;al.$ \cite{PGF13} studied that the class
$\Omega$ enlarges the class $\Phi$, $i.e.$, $\Phi\subset \Omega$.\\

\noindent\textbf{Definition 21}. Let $X$ be a nonempty set and $f$
and $g$ two self-mappings on $X$. Then
 an element $x\in X$ is called a
coincidence point of $f$ and $g$ if
$$f(x)=g(x)=\overline{x},$$
for some $\overline{x}\in X$. Moreover, $\overline{x}$ is called a
point of coincidence of $f$ and $g$. Furthermore, if
$\overline{x}=x$, then $x$ is called a
common fixed point of $f$ and $g$.\\

The following coincidence theorems are crucial results to prove our main results.\\

\noindent\textbf{Lemma 1.} Let $(X,d,\preceq)$ be an ordered metric
space and $E$ an $\overline{\rm O}$-complete (resp. $\underline{\rm
O}$-complete) subspace of $X$. Let $f$ and $g$ be two self-mappings
on $X$. Suppose that the following conditions hold:
\begin{enumerate}
\item [{(i)}] $f(X)\subseteq g(X)\cap E$,
\item [{(ii)}] $f$ is $g$-increasing,
\item [{(iii)}] $f$ and $g$ are $\overline{\rm O}$-compatible (resp. $\underline{\rm O}$-compatible),
\item [{(iv)}] $g$ is $\overline{\rm O}$-continuous (resp. $\underline{\rm O}$-continuous),
\item [{(v)}] either $f$ is $\overline{\rm O}$-continuous (resp. $\underline{\rm O}$-continuous) or $(E,d,\preceq)$ has {\it g-ICU} property (resp. {\it g-DCL} property),
\item [{(vi)}] there exists $x_{0}\in X$ such that
$g(x_{0})\preceq f(x_{0})$ (resp. $g(x_{0})\succeq f(x_{0})$),
\item [{(vii)}] there exists $\varphi\in \Omega$ such that
$$d(fx,fy)\leq\varphi(d(gx,gy))\;\;\forall~x,y\in X~{\rm with}~
g(x)\prec\succ g(y).$$
 \end{enumerate}
Then $f$ and $g$ have a coincidence point. Moreover, if the
following condition also holds:\\
\indent\hspace{5mm}(viii) for each pair $x,y\in X$, $\exists~z\in X$ such that $g(x)\prec\succ g(z)$ and $g(y)\prec\succ g(z)$,\\
then $f$ and $g$ have a unique point of coincidence, which remains also a unique common fixed point.\\

\noindent\textbf{Lemma 2.} Let $(X,d,\preceq)$ be an ordered metric
space and $E$ an $\overline{\rm O}$-complete (resp. $\underline{\rm
O}$-complete) subspace of $X$. Let $f$ and $g$ be two self-mappings
on $X$. Suppose that the following conditions hold:
\begin{enumerate}
\item [{(i)}] $f(X)\subseteq E\subseteq g(X)$,
\item [{(ii)}] $f$ is $g$-increasing,
\item [{(iii)}] either $f$ is $(g,\overline{\rm O})$-continuous (resp. $(g,\underline{\rm O})$-continuous) or $f$ and $g$ are continuous or $(E,d,\preceq)$ has
{\it ICU} property (resp. {\it DCL} property),
\item [{(iv)}] there exists $x_{0}\in X$ such that
$g(x_{0})\preceq f(x_{0})$ (resp. $g(x_{0})\succeq f(x_{0})$),
\item [{(v)}] there exists $\varphi\in \Omega$ such that
$$d(fx,fy)\leq\varphi(d(gx,gy))\;\;\forall~x,y\in X~{\rm with}~
g(x)\prec\succ g(y).$$
 \end{enumerate}
Then $f$ and $g$ have a coincidence point. Moreover, if the
following condition also holds:\\
\indent\hspace{5mm}(vi) for each pair $x,y\in X$, $\exists~z\in X$ such that $g(x)\prec\succ g(z)$ and $g(y)\prec\succ g(z)$,\\
then $f$ and $g$ have a unique point of coincidence.\\

We skip the proofs of above lemmas as they are proved in Alam $et\;al.$ \cite{PGF13,PGF17,PGF18}.

\section{Extended notions upto product sets}
Recall that a binary operation $\ast$ on a set $S$ is a mapping from
$S\times S$ to $S$ and a permutation $\pi$ on a set $S$ is a one-one
mapping from a $S$ onto itself ($cf.$ Herstein \cite{ALG}).
Throughout this manuscript, we adopt the following notations:
\begin{enumerate}
\item [{(1)}] In order to understand a binary operation $\ast$ on $I_n$, we denote the image of any element $(i,k)\in I_n\times I_n$ under $\ast$ by $i_k$ rather
than $\ast(i,k)$.
\item [{(2)}] A binary operation $\ast$ on $I_n$ can be identically
represented by an $n\times n$ matrix throughout its ordered image
such that the first and second components run over rows and columns
respectively, $i.e.$,
$$\ast=[m_{ik}]_{n\times n}\; {\rm where}\; m_{ik}=i_k \;{\rm for\;each}\; i,k\in I_n.$$
\item [{(3)}] A permutation $\pi$ on $I_n$ can be identically represented by an $n$-tuple throughout its ordered image,
$i.e.$,
$$\pi=(\pi(1),\pi(2),...,\pi(n)).$$
\item [{(4)}] $\mathfrak{B_{n}}$ denotes the family of all binary operations $\ast$ on
$I_n$, $i.e.,$
$$\mathfrak{B_{n}}=\{\ast:\ast:I_n\times I_n\rightarrow I_n\}.$$
\end{enumerate}
\begin{rem} It is clear for each $i\in I_n$ that
$$\{i_1,i_2,...,i_n\}\subseteq I_n.$$
\end{rem}

\noindent\textbf{Definition 19} (Alam $et\;al.$ \cite{unif1}). Let
$X$ be a nonempty set, $\ast\in\mathfrak{B_{n}}$ and $F:X^{n}
\rightarrow X$ and $g: X\rightarrow X$ two mappings. An element
$(x_{1},x_{2},...,x_{n})\in X^{n}$  is called an $n$-tupled
coincidence point of $F$ and $g$ w.r.t. $\ast$ (or, in short,
$\ast$-coincidence point of $F$ and $g$) if
$$F(x_{i_1},x_{i_2},...,x_{i_n})=g(x_i)\;\;{\rm for\;each}\;i\in I_n.$$
In this case $(gx_1,gx_2,...,gx_n)$ is called point of
$\ast$-coincidence of $F$ and $g$.\\
Notice that if $g$ is an identity mapping on $I_n$ then the notion
employed in Definition 19 is called an $n$-tupled fixed point of $F$
w.r.t. $\ast$ (or, in short,
$\ast$-fixed point of $F$).\\

\noindent\textbf{Definition 20} (Alam $et\;al.$ \cite{unif1}). Let
$X$ be a nonempty set, $\ast\in\mathfrak{B_{n}}$ and $F:X^{n}
\rightarrow X$ and $g: X\rightarrow X$ two mappings. An element
$(x_{1},x_{2},...,x_{n})\in X^{n}$  is called a common $n$-tupled
fixed point of $F$ and $g$ w.r.t. $\ast$ (or, in short, common
$\ast$-fixed point of $F$ and $g$) if
$$F(x_{i_1},x_{i_2},...,x_{i_n})=g(x_i)=x_i\;\;{\rm for\;each}\;i\in I_n.$$

In the following lines, we define four special types $n$-tupled
fixed points, which are somewhat natural.\\

\noindent\textbf{Definition 21} (Alam $et\;al.$ \cite{unif1}). Let
$X$ be a nonempty set and $F:X^{n} \rightarrow X$ a mapping. An
element $(x_{1},x_{2},...,x_{n})\in X^{n}$ is called a forward
cyclic $n$-tupled fixed point of $F$ if
$$F(x_{i},x_{i+1},...,x_n,x_1,...,x_{i-1})=x_i\;{\rm for\;each}\;i\in I_n$$
$i.e.$\\
\indent\hspace{5cm}$F(x_{1},x_{2},...,x_{n})=x_{1},$\\
\indent\hspace{5cm}$F(x_{2},x_{3},...,x_{n},x_{1})=x_{2},$\\
\indent\hspace{5.6cm}$\vdots$\\
\indent\hspace{5cm}$F(x_{n},x_{1},x_{2},...,x_{n-1})=x_{n}.$\\
This was initiated by Samet and Vetro \cite{NFI}. To obtain this we define $\ast$ as\\
\vspace{0.2cm}\\
\indent\hspace{3cm}$i_k={\begin{cases}i+k-1\;\;\;\;\;\;\;\;\;\;1\leq
k\leq{n-i+1}\cr
\hspace{0.0in}i+k-n-1\;\;\;\;{n-i+2}\leq k\leq{n}\cr\end{cases}}$\\
\vspace{0.2cm}\\
$i.e.$\indent\hspace{2.6cm}$\ast=\left[\begin{matrix}
1 &2 &\cdots &{n-1} &n\\
2 &3 &\cdots &{n} &{1}\\
\vdots &\vdots &&\vdots &\vdots \\
n &1 &\cdots &{n-2} &{n-1}\\
\end{matrix}\right]_{n\times n}$\\

\noindent\textbf{Definition 22} (Alam $et\;al.$ \cite{unif1}). Let
$X$ be a nonempty set and $F:X^{n} \rightarrow X$ a mapping. An
element $(x_{1},x_{2},...,x_{n})\in X^{n}$  is called a backward
cyclic $n$-tupled fixed point of $F$ if
$$F(x_{i},x_{i-1},...,x_{1},x_n,x_{n-1},...,x_{i+1})=x_i\;{\rm for\;each}\;i\in I_n$$
$i.e.$\\
\indent\hspace{5cm}$F(x_{1},x_{n},x_{n-1},...,x_{2})=x_{1},$\\
\indent\hspace{5cm}$F(x_{2},x_{1},x_{n},...,x_{3})=x_{2},$\\
\indent\hspace{5.6cm}$\vdots$\\
\indent\hspace{5cm}$F(x_{n},x_{n-1},x_{n-2},...,x_{1})=x_{n}.$\\
To obtain this we define $\ast$ as\\
\vspace{0.2cm}\\
\indent\hspace{3cm}$i_k={\begin{cases}i-k+1\;\;\;\;\;\;\;\;\;\;1\leq
k\leq{i}\cr
\hspace{0.0in}n+i-k+1\;\;\;\;{i+1}\leq k\leq{n-1}\cr\end{cases}}$\\
\vspace{0.2cm}\\
$i.e.$\indent\hspace{2.6cm}$\ast=\left[\begin{matrix}
1 &n &{n-1} &\cdots &2\\
2 &1 &{n} &\cdots &3\\
\vdots &\vdots &\vdots &&\vdots \\
n &{n-1} &{n-2} &\cdots &1\\
\end{matrix}\right]_{n\times n}$\\

\noindent\textbf{Definition 23} (Alam $et\;al.$ \cite{unif1}). Let
$X$ be a nonempty set and $F:X^{n} \rightarrow X$ a mapping. An
element $(x_{1},x_{2},...,x_{n})\in X^{n}$  is called a 1-skew
cyclic $n$-tupled fixed point of $F$ if
$$F(x_{i},x_{i-1},...,x_2,x_1,x_2,...,x_{n-i+1})=x_i\;{\rm for\;each}\;i\in I_n.$$
This was introduced by Gordji and Ramezani \cite{NX1}. To find this we define $\ast$ as\\
\vspace{0.2cm}\\
\indent\hspace{3cm}$i_k={\begin{cases}i-k+1\;\;\;\;\;\;\;\;\;\;1\leq
k\leq{i}\cr
\hspace{0.0in}k-i+1\;\;\;\;{i+1}\leq k\leq{n}\cr\end{cases}}$\\

\noindent\textbf{Definition 24} (Alam $et\;al.$ \cite{unif1}). Let
$X$ be a nonempty set and $F:X^{n} \rightarrow X$ a mapping. An
element $(x_{1},x_{2},...,x_{n})\in X^{n}$ is called a $n$-skew
cyclic $n$-tupled fixed point of $F$ if
$$F(x_{i},x_{i+1},...,x_{n-1},x_n,x_{n-1},...,x_{n-i+1})=x_i\;{\rm for\;each}\;i\in I_n.$$
To find this we define $\ast$ as\\
\indent\hspace{3cm}$i_k={\begin{cases}i+k-1\;\;\;\;\;\;\;\;\;\;1\leq
k\leq{n-i+1}\cr \hspace{0.0in}2n-i-k+1\;\;\;\;{n-i+2}\leq
k\leq{n}\cr\end{cases}}$

\noindent\textbf{Definition 25} (Alam $et\;al.$ \cite{unif1}). A
binary operation $\ast$ on $I_n$ is called permuted if each row of
matrix representation of $\ast$ forms a permutation on $I_n.$
\begin{example} (Alam $et\;al.$ \cite{unif1}). On $I_3$, consider two binary operations\\
\indent\hspace{2cm}$\ast$=$\left[\begin{matrix}
1 &2 &3\\
2 &1 &3\\
3 &2 &1\\
\end{matrix}\right]$,\;\;
$\circ$=$\left[\begin{matrix}
1 &2 &3\\
2 &1 &3\\
3 &3 &2\\
\end{matrix}\right]$\\
$\ast$ is permuted as each of rows $(1,2,3),(2,1,3),(3,2,1)$ is a
permutation on $I_3$. While $\circ$ is not permuted as last row
$(3,3,2)$ is not permutation on $I_3$.
\end{example}
It is clear that binary operations defined for forward cyclic
and backward cyclic $n$-tupled fixed points are permuted while
for 1-skew cyclic and $n$-skew cyclic $n$-tupled fixed
points are not permuted.
\begin{prop} (Alam $et\;al.$ \cite{unif1}). A permutation $\ast$ on $I_n$ is permuted iff for each $i\in I_n,$
$$\{i_1,i_2,...,i_n\}=I_n.$$
\end{prop}

\noindent\textbf{Definition 28} (Alam $et\;al.$ \cite{unif1}). Let $(X,d)$ be a metric space,
$F:X^{n} \rightarrow X$ a mapping and $(x_1,x_2,...,x_n)\in X^n$. We
say that $F$ is continuous at $(x_1,x_2,...,x_n)$ if for any
sequences $\{x_1^{(m)}\},\{x_2^{(m)}\},$ $..., \{x_n^{(m)}\}\subset
X$,
$$x_1^{(m)}\stackrel{d}{\longrightarrow} x_1,\;x_2^{(m)}\stackrel{d}{\longrightarrow} x_2,...,\;x_n^{(m)}\stackrel{d}{\longrightarrow} x_n$$
$$\Longrightarrow F(x_1^{(m)},x_2^{(m)},...,x_n^{(m)})\stackrel{d}{\longrightarrow}F(x_1,x_2,...,x_n).$$
Moreover, $F$ is called continuous if it is continuous at each point of $X^n$.\\

\noindent\textbf{Definition 29} (Alam $et\;al.$ \cite{unif1}). Let $(X,d)$ be a metric space and
$F:X^{n} \rightarrow X$ and $g: X\rightarrow X$ two mappings and
$(x_1,x_2,...,x_n)\in X^n$. We say that $F$ is $g$-continuous at
$(x_1,x_2,...,x_n)$ if for any sequences
$\{x_1^{(m)}\},\{x_2^{(m)}\},...,\{x_n^{(m)}\}\subset X$,
$$g(x_1^{(m)})\stackrel{d}{\longrightarrow} g(x_1),\;g(x_2^{(m)})\stackrel{d}{\longrightarrow} g(x_2),...,\;g(x_n^{(m)})\stackrel{d}{\longrightarrow} g(x_n)$$
$$\Longrightarrow F(x_1^{(m)},x_2^{(m)},...,x_n^{(m)})\stackrel{d}{\longrightarrow}F(x_1,x_2,...,x_n).$$
Moreover, $F$ is called $g$-continuous if it is $g$-continuous at each point of $X^n$.\\

Notice that setting $g=I$ (identity mapping on $X$), Definition 29
reduces to Definition 28.\\

\noindent\textbf{Definition 30}. Let $(X,d,\preceq)$ be an ordered
metric space, $F:X^{n} \rightarrow X$ a mapping and
$(x_1,x_2,...,x_n)\in X^n$. We say that $F$ is
\begin{enumerate}
\item[{(i)}] $\overline{\rm O}$-continuous at $(x_1,x_2,...,x_n)$ if for any sequences
$\{x_1^{(m)}\},\{x_2^{(m)}\},...,\{x_n^{(m)}\}\subset X$,
$$x_1^{(m)}\uparrow x_1,\;x_2^{(m)}\uparrow x_2,...,\;x_n^{(m)}\uparrow x_n$$
$$\Longrightarrow F(x_1^{(m)},x_2^{(m)},...,x_n^{(m)})\stackrel{d}{\longrightarrow}F(x_1,x_2,...,x_n),$$
\item[{(ii)}] $\underline{\rm O}$-continuous at $(x_1,x_2,...,x_n)$ if for any sequences
$\{x_1^{(m)}\},\{x_2^{(m)}\},...,\{x_n^{(m)}\}\subset X$,
$$x_1^{(m)}\downarrow x_1,\;x_2^{(m)}\downarrow x_2,...,\;x_n^{(m)}\downarrow x_n$$
$$\Longrightarrow F(x_1^{(m)},x_2^{(m)},...,x_n^{(m)})\stackrel{d}{\longrightarrow}F(x_1,x_2,...,x_n)\;{\rm and}$$
\item[{(iii)}] O-continuous at
$(x_1,x_2,...,x_n)$ if for any sequences
$\{x_1^{(m)}\},\{x_2^{(m)}\},...,\{x_n^{(m)}\}\subset X$,
$$x_1^{(m)}\uparrow\downarrow x_1,\;x_2^{(m)}\uparrow\downarrow x_2,...,\;x_n^{(m)}\uparrow\downarrow x_n$$
$$\Longrightarrow F(x_1^{(m)},x_2^{(m)},...,x_n^{(m)})\stackrel{d}{\longrightarrow}F(x_1,x_2,...,x_n).$$
\end{enumerate}
Moreover, $F$ is called O-continuous (resp. $\overline{\rm
O}$-continuous, $\underline{\rm O}$-continuous) if it is
O-continuous (resp. $\overline{\rm O}$-continuous, $\underline{\rm
O}$-continuous) at each point of $X^n$.
\begin{rem} In an ordered
metric space, continuity $\Rightarrow$ O-continuity $\Rightarrow$
$\overline{\rm O}$-continuity as well as $\underline{\rm
O}$-continuity .\end{rem}

\noindent\textbf{Definition 31}. Let $(X,d,\preceq)$ be an ordered
metric space, $F:X^{n} \rightarrow X$ and $g: X\rightarrow X$ two
mappings and $(x_1,x_2,...,x_n)\in X^n$. We say that $F$ is
\begin{enumerate}
\item[{(i)}] $(g,\overline{\rm O})$-continuous at $(x_1,x_2,...,x_n)$ if for any sequences
$\{x_1^{(m)}\},\{x_2^{(m)}\},...,\{x_n^{(m)}\}\subset X$,
$$g(x_1^{(m)})\uparrow g(x_1),\;g(x_2^{(m)})\uparrow g(x_2),...,\;g(x_n^{(m)})\uparrow g(x_n)$$
$$\Longrightarrow F(x_1^{(m)},x_2^{(m)},...,x_n^{(m)})\stackrel{d}{\longrightarrow}F(x_1,x_2,...,x_n),$$
\item[{(ii)}] $(g,\underline{\rm O})$-continuous at $(x_1,x_2,...,x_n)$ if for any sequences
$\{x_1^{(m)}\},\{x_2^{(m)}\},...,\{x_n^{(m)}\}\subset X$,
$$g(x_1^{(m)})\downarrow g(x_1),\;g(x_2^{(m)})\downarrow g(x_2),...,\;g(x_n^{(m)})\downarrow g(x_n)$$
$$\Longrightarrow F(x_1^{(m)},x_2^{(m)},...,x_n^{(m)})\stackrel{d}{\longrightarrow}F(x_1,x_2,...,x_n)\;{\rm and}$$
\item[{(iii)}] $(g,{\rm O})$-continuous at $(x_1,x_2,...,x_n)$ if for any
sequences $\{x_1^{(m)}\},\{x_2^{(m)}\},...,\{x_n^{(m)}\}\subset X$,
$$g(x_1^{(m)})\uparrow\downarrow g(x_1),\;g(x_2^{(m)})\uparrow\downarrow g(x_2),...,\;g(x_n^{(m)})\uparrow\downarrow g(x_n)$$$$\Longrightarrow
F(x_1^{(m)},x_2^{(m)},...,x_n^{(m)})\stackrel{d}{\longrightarrow}F(x_1,x_2,...,x_n).$$
\end{enumerate}

Notice that setting $g=I$ (identity mapping on $X$), Definition 31
reduces to Definition 30.
\begin{rem} In an ordered
metric space, $g$-continuity $\Rightarrow$ $(g,{\rm O})$-continuity
$\Rightarrow$ $(g,\overline{\rm O})$-continuity as well as
$(g,\underline{\rm O})$-continuity .\end{rem}

\noindent\textbf{Definition 32} (Alam $et\;al.$ \cite{unif1}). Let $X$ be a nonempty set and
$F:X^{n} \rightarrow X$ and $g: X\rightarrow X$ two mappings. We say
that $F$ and $g$ are commuting if for all $x_1,x_2,...,x_n\in X,$
$$g(F(x_1,x_2,...,x_n))=F(gx_1,gx_2,...,gx_n).$$

\noindent\textbf{Definition 33} (Alam $et\;al.$ \cite{unif1}). Let $(X,d)$ be a metric space and
$F:X^{n} \rightarrow X$ and $g: X\rightarrow X$ two mappings. We say
that $F$ and $g$ are $\ast$-compatible if for any sequences
$\{x_1^{(m)}\},\{x_2^{(m)}\},...,\{x_n^{(m)}\}\subset X$ and for any
$z_1,z_2,...,z_n\in X$,
$$g(x_i^{(m)})\stackrel{d}{\longrightarrow} z_i\;{\rm and}\;F(x_{i_1}^{(m)},x_{i_2}^{(m)},...,x_{i_n}^{(m)})\stackrel{d}{\longrightarrow} z_i\;\;{\rm for~each}~i\in I_n$$
$$\Longrightarrow\lim\limits_{m\to \infty}d(gF(x_{i_1}^{(m)},x_{i_2}^{(m)},...,x_{i_n}^{(m)}),F(gx_{i_1}^{(m)},gx_{i_2}^{(m)},...,gx_{i_n}^{(m)}))=0\;\;{\rm for~each}~i\in I_n.$$

\noindent\textbf{Definition 34}. Let $(X,d,\preceq)$ be an ordered
metric space and $F:X^{n} \rightarrow X$ and $g: X\rightarrow X$ two
mappings. We say that $F$ and $g$ are
\begin{enumerate}
\item[{(i)}] $(\ast,\overline{\rm O})$-compatible if for any sequences
$\{x_1^{(m)}\},\{x_2^{(m)}\},...,\{x_n^{(m)}\}\subset X$ and for any
$z_1,z_2,...,z_n\in X$,
$$g(x_i^{(m)})\uparrow z_i\;{\rm and}\;F(x_{i_1}^{(m)},x_{i_2}^{(m)},...,x_{i_n}^{(m)})\uparrow z_i\;\;{\rm for~each}~i\in I_n$$
$$\Longrightarrow\lim\limits_{m\to \infty}d(gF(x_{i_1}^{(m)},x_{i_2}^{(m)},...,x_{i_n}^{(m)}),F(gx_{i_1}^{(m)},gx_{i_2}^{(m)},...,gx_{i_n}^{(m)}))=0,$$
\item[{(ii)}] $(\ast,\underline{\rm O})$-compatible if for any sequences
$\{x_1^{(m)}\},\{x_2^{(m)}\},...,\{x_n^{(m)}\}\subset X$ and for any
$z_1,z_2,...,z_n\in X$,
$$g(x_i^{(m)})\downarrow z_i\;{\rm and}\;F(x_{i_1}^{(m)},x_{i_2}^{(m)},...,x_{i_n}^{(m)})\downarrow z_i\;\;{\rm for~each}~i\in I_n$$
$$\Longrightarrow\lim\limits_{m\to \infty}d(gF(x_{i_1}^{(m)},x_{i_2}^{(m)},...,x_{i_n}^{(m)}),F(gx_{i_1}^{(m)},gx_{i_2}^{(m)},...,gx_{i_n}^{(m)}))=0\;\;{\rm for~each}~i\in I_n\;{\rm and}$$
\item[{(iii)}] $(\ast,{\rm O})$-compatible if for any sequences
$\{x_1^{(m)}\},\{x_2^{(m)}\},...,\{x_n^{(m)}\}\subset X$ and for any
$z_1,z_2,...,z_n\in X$,
$$g(x_i^{(m)})\uparrow\downarrow z_i\;{\rm and}\;F(x_{i_1}^{(m)},x_{i_2}^{(m)},...,x_{i_n}^{(m)})\uparrow\downarrow z_i\;\;{\rm for~each}~i\in I_n$$
$$\Longrightarrow\lim\limits_{m\to \infty}d(gF(x_{i_1}^{(m)},x_{i_2}^{(m)},...,x_{i_n}^{(m)}),F(gx_{i_1}^{(m)},gx_{i_2}^{(m)},...,gx_{i_n}^{(m)}))=0\;\;{\rm for~each}~i\in I_n.$$

\end{enumerate}

\noindent\textbf{Definition 35} (Alam $et\;al.$ \cite{unif1}). Let
$X$ be a nonempty set and $F:X^{n} \rightarrow X$ and $g:
X\rightarrow X$ two mappings. We say that $F$ and $g$ are
$(\ast,w)$-compatible if for any $x_1,x_2,...,x_n\in X,$
$$g(x_i)=F(x_{i_1},x_{i_2},...,x_{i_n})\;\;{\rm for~each}~i\in I_n$$
$$\Longrightarrow g(F(x_{i_1},x_{i_2},...,x_{i_n}))=F(gx_{i_1},gx_{i_2},...,gx_{i_n})\;\;{\rm for~each}~i\in I_n.$$
\begin{rem} Evidently, in an ordered metric space, commutativity $\Rightarrow$
$\ast$-compatibility $\Rightarrow$ $(\ast,{\rm O})$-compatibility
$\Rightarrow$ $(\ast,\overline{\rm O})$-compatibility as well as
$(\ast,\underline{\rm O})$-compatibility $\Rightarrow$
$(\ast,w)$-compatibility for a pair of mappings $F:X^{n} \rightarrow
X$ and $g: X\rightarrow X$.
\end{rem}
\begin{prop} (Alam $et\;al.$ \cite{unif1}). If $F$ and $g$ are
$(\ast,w)$-compatible, then every point of $\ast$-coincidence of $F$
and $g$ is also an $\ast$-coincidence point of $F$
 and $g$.\end{prop}

\section{Auxiliary results}

In this section, we discuss some basic results, which provide the
tools for reduction of the multi-tupled fixed point results from the
corresponding fixed point results. Before doing this, we consider
the following induced notations:
\begin{enumerate}
\item [{(1)}] For any U$=(x_{1},x_{2},...,x_{n})\in X^n$, for an $\ast\in \mathfrak{B_{n}}$ and for each $i\in I_n$, U$^{\ast}_i$ denotes the ordered element
 $(x_{i_1},x_{i_2},...,x_{i_n})$ of $X^n$.
\item [{(2)}] For each $\ast\in \mathfrak{B_{n}}$, a mapping $F:X^n\rightarrow X$ induce an associated mapping $F_\ast:X^n\rightarrow X^n$ defined by
$$F_\ast({\rm U})=(F{\rm U}^{\ast}_1,F{\rm U}^{\ast}_2,...,F{\rm U}^{\ast}_n)\;\;\forall~{\rm U}\in X^n.$$
\item [{(3)}] A mapping $g:X\rightarrow X$ induces an associated mapping $G:X^n\rightarrow
X^n$ defined by
$$G({\rm U})=(gx_1,gx_2,...,gx_n)\;\;\forall~{\rm U}=(x_{1},x_{2},...,x_{n})\in X^n. $$
\item [{(4)}] For a metric space $(X,d)$, $\Delta_n$ and $\nabla_n$ denote
two metrics on product set $X^n$ defined by:\\ for all
U=$(x_{1},x_{2},...,x_{n})$, V=$(y_{1},y_{2},...,y_{n})\in X^n,$
$$\Delta_n({\rm U,V})=\frac{1}{n}\sum\limits_{i=1}^{n}d(x_i,y_i)$$
$$\nabla_n({\rm U,V})=\max\limits_{i\in I_n}d(x_i,y_i).$$
\item [{(5)}] For any ordered set $(X,\preceq)$, $\sqsubseteq_n$ denotes a partial order on $X^n$ defined
by:\\
for all U=$(x_{1},x_{2},...,x_{n})$, V=$(y_{1},y_{2},...,y_{n})\in
X^n,$
$${\rm U}\sqsubseteq_n{\rm V}\Leftrightarrow x_i\preceq y_i \;{\rm {for\; each} }\; i\in I_n.$$
\end{enumerate}
\begin{rem} The following facts are straightforward:
\begin{enumerate}
\item [{(i)}] $F_\ast(X^n)\subseteq (FX^n)^n.$
\item [{(ii)}] $G(X^n)=(gX)^n.$
\item [{(iii)}] $(G{\rm U})^\ast_i=G({\rm U}^\ast_i)\;\forall~{\rm U}\in X^n.$
\item [{(iv)}] $\frac{1}{n}\nabla_n\leq\Delta_n\leq\nabla_n$ ($i.e.$ both the
metrics $\Delta_n$ and $\nabla_n$ are equivalent).
\end{enumerate}
\end{rem}
\noindent{\bf Lemma 3} (Alam $et\;al.$ \cite{unif1}). Let $X$ be a nonempty set, $E\subseteq X$,
$F:X^{2} \rightarrow X$ and $g: X\rightarrow X$ two mappings and
$\ast\in\mathfrak{B_{n}}$.
\begin{enumerate}
\item [{(i)}] If $F(X^n)\subseteq g(X)\cap E$ then $F_\ast(X^n)\subseteq (FX^n)^n\subseteq
G(X^n)\cap E^n$.
\item [{(ii)}] If $F(X^n)\subseteq E\subseteq g(X)$ then $F_\ast(X^n)\subseteq
(FX^n)^n\subseteq E^n \subseteq G(X^n)$.
\item [{(iii)}] An element $(x_{1},x_{2},...,x_{n})\in
X^n$ is $\ast$-coincidence point of $F$ and $g$ iff
$(x_{1},x_{2},...,x_{n})$ is a coincidence point of $F_{\ast}$ and $G$.
\item [{(iv)}] An element $(\overline{x}_{1},\overline{x}_{2},...,\overline{x}_{n})\in
X^n$ is point of $\ast$-coincidence of $F$ and $g$ iff
$(\overline{x}_{1},\overline{x}_{2},...,\overline{x}_{n})$ is a
point of coincidence of $F_{\ast}$ and $G$.
\item [{(v)}] An element $(x_{1},x_{2},...,x_{n})\in X^n$ is common $\ast$-fixed point of $F$ and $g$ iff
$(x_{1},x_{2},$$...,$$x_{n})$ is a common fixed point of $F_{\ast}$ and $G$.
\end{enumerate}

\noindent{\bf Lemma 4.} Let $(X,\preceq)$ be an ordered set, $g:
X\rightarrow X$ a mapping and $\ast\in\mathfrak{B_{n}}$. If $G({\rm
U})\sqsubseteq_n G({\rm V})$ for some U,V$\in X^n$ then for each
$i\in I_n,$ $G({\rm U}^\ast_i)\sqsubseteq_n G({\rm
V}^\ast_i)$.\\

\noindent{\bf Proof.} Let U=$(x_{1},x_{2},...,x_{n})$ and
V=$(y_{1},y_{2},...,y_{n})$ be such that $G({\rm U})\sqsubseteq_n
G({\rm V})$, then we have
$$(gx_{1},gx_{2},...,gx_{n})\sqsubseteq_{n}(gy_{1},gy_{2},...,gy_{n}),$$
$$\Rightarrow g(x_i)\preceq g(y_i) \;{\rm {for\; each} }\; i\in I_n$$
$$\Rightarrow g(x_{i_k})\preceq g(y_{i_k})\;{\rm {for\; each} }\;i\in I_n{\rm {and~for\; each} }\;k\in I_n$$
$$\Rightarrow (gx_{i_1},gx_{i_2},...,gx_{i_n})\sqsubseteq_{n}(gy_{i_1},gy_{i_2},...,gy_{i_n})\;\;{\rm for~ each~} i\in I_n,$$
$i.e.$ $$G({\rm U}^\ast_i)\sqsubseteq_{n}G({\rm V}^\ast_i)\;\;{\rm
for~ each~} i\in I_n.$$

\noindent{\bf Lemma 5.} Let $(X,\preceq)$ be an ordered set,
$F:X^{2} \rightarrow X$ and $g: X\rightarrow X$ two mappings and
$\ast\in\mathfrak{B_{n}}$. If $F$ has
$g$-monotone property then $F_{\ast}$ is $G$-increasing in ordered set $(X^n,\sqsubseteq_{n})$.\\
\noindent{\bf Proof.} Take U=$(x_{1},x_{2},...,x_{n})$,
V=$(y_{1},y_{2},...,y_{n})\in X^n$ with $G({\rm
U})\sqsubseteq_{n}G({\rm V})$. Using Lemma 4, we obtain
$$G({\rm U}^\ast_i)\sqsubseteq_n G({\rm
V}^\ast_i) \;{\rm {for\; each} }\; i\in I_n,$$
which implies, for all $i\in I_n$, that
$$\Rightarrow g(x_{i_k})\preceq g(y_{i_k})\;{\rm {for\; each}}\;k\in I_n\eqno(1)$$
On using (1) and  $g$-monotone property of $F$, we obtain, for all
$i\in I_n$, that
$$F(x_{i_1},x_{i_2},...,x_{i_n})\preceq F(y_{i_1},y_{i_2},...,y_{i_n})$$
$i.e.$
$$F({\rm U}^\ast_i)\preceq F({\rm V}^\ast_i).\eqno(2)$$
Using (2), we get
\begin{eqnarray*}
F_\ast({\rm U})&=&(F{\rm U}^{\ast}_1,F{\rm U}^{\ast}_2,...,F{\rm U}^{\ast}_n)\\
&\sqsubseteq_{n}& (F{\rm V}^{\ast}_1,F{\rm V}^{\ast}_2,...,F{\rm V}^{\ast}_n)\\
&=& F_\ast({\rm V}).
\end{eqnarray*}
Hence, $F_\ast$ is $G$-increasing.\\

\noindent{\bf Lemma 6} (Alam $et\;al.$ \cite{unif1}). Let $(X,d)$ be
a metric space, $g: X\rightarrow X$ a mapping and
$\ast\in\mathfrak{B_{n}}$. Then, for any
U=$(x_{1},x_{2},...,x_{n})$,V=$(y_{1},y_{2},...,y_{n})\in X^n$ and
for each $i\in I_n$,
\begin{enumerate}
\item [{(i)}] $\frac{1}{n}\sum\limits_{k=1}^{n}d(gx_{i_k},gy_{i_k})=\frac{1}{n}\sum\limits_{j=1}^{n}d(gx_j,gy_j)=\Delta_n(G{\rm U},G{\rm V})$\; provided $\ast$
is permuted,\\
\item [{(ii)}] $\max\limits_{k\in
I_n}d(gx_{i_k},gy_{i_k})=\max\limits_{j\in
I_n}d(gx_j,gy_j)=\nabla_n(G{\rm U},G{\rm V})$\;provided $\ast$
is permuted,\\
\item [{(iii)}] $\max\limits_{k\in
I_n}d(gx_{i_k},gy_{i_k})\leq\max\limits_{j\in
I_n}d(gx_j,gy_j)=\nabla_n(G{\rm U},G{\rm V})$.
\end{enumerate}

\begin{prop} (Alam $et\;al.$ \cite{unif1}). Let $(X,d)$ be a metric space. Then for any sequence ${\ \rm U}^{(m)}\ \subset X^n$
and any ${\rm U}\in X^n$, where
${\rm U}^{(m)}=(x^{(m)}_1,x^{(m)}_2,...,x^{(m)}_n)$ and
${\rm U}=(x_1,x_2,...,x_n)$
\begin{enumerate}
\item [{(i)}] ${\rm U}^{(m)}\stackrel{\Delta_n}{\longrightarrow} {\rm U}\Leftrightarrow x_i^{(m)}\stackrel{d}{\longrightarrow} x_i\;{\rm for~each}~i\in I_n.$
\item [{(ii)}] ${\rm U}^{(m)}\stackrel{\nabla_n}{\longrightarrow} {\rm U}\Leftrightarrow x_i^{(m)}\stackrel{d}{\longrightarrow} x_i\;{\rm for~each}~i\in I_n.$
\end{enumerate}
\end{prop}

\noindent{\bf Lemma 7} (Alam $et\;al.$ \cite{unif1}). Let $(X,d)$ be
a metric space, $F:X^{n} \rightarrow X$ and $g:X\rightarrow X$ two
mappings and $\ast\in\mathfrak{B_{n}}$.
\begin{enumerate}
\item [{(i)}] If $g$ is continuous then $G$ is continuous in both
metric spaces $(X^n,\Delta_n)$ and $(X^n,\nabla_n)$,
\item [{(ii)}] If $F$ is continuous then $F_{\ast}$ is continuous in
both metric spaces $(X^n,\Delta_n)$ and $(X^n,\nabla_n)$.
\end{enumerate}

\begin{prop} (Alam $et\;al.$ \cite{unif1}) Let $(X,d,\preceq)$ be an ordered metric space and $\{{\rm
U}^{(m)}\}$ a sequence in $X^n,$ where ${\rm
U}^{(m)}=(x_1^{(m)},x_2^{(m)},...,x_n^{(m)})$.
\begin{enumerate}
\item [{(i)}] If $\{{\rm
U}^{(m)}\}$ is increasing (resp. decreasing) in
$(X^n,\sqsubseteq_{n})$ then each of
$\{x_1^{(m)}\}$,$\{x_2^{(m)}\}$,...,$\{x_n^{(m)}\}$ is increasing
(resp. decreasing) in $(X,\preceq)$.
\item [{(ii)}] If $\{{\rm
U}^{(m)}\}$ is Cauchy in $(X^n,\Delta_n)$ (similarly in
$(X^n,\nabla_n)$) then each of $\{x_1^{(m)}\}$,$\{x_2^{(m)}\}$,
...,$\{x_n^{(m)}\}$ is Cauchy in $(X,d)$.
\end{enumerate}
\end{prop}
{\noindent\bf{Lemma 8.}} Let $(X,d,\preceq)$ be an ordered metric
space, $E\subseteq X$, $F:X^{n} \rightarrow X$ and $g:X\rightarrow
X$ two mappings and $\ast\in\mathfrak{B_{n}}$.
\begin{enumerate}
\item [{(i)}] If $(E,d,\preceq)$ is $\overline{\rm O}$-complete (resp. $\underline{\rm O}$-complete) then
$(E^n,\Delta_n,\sqsubseteq_{n})$ and
$(E^n,\nabla_n,\sqsubseteq_{n})$ both are $\overline{\rm
O}$-complete (resp. $\underline{\rm O}$-complete).
\item [{(ii)}] If $F$ and $g$ are $(\ast,\overline{\rm
O})$-compatible pair (resp. $(\ast,\underline{\rm O})$-compatible
pair) then $F_{\ast}$ and $G$ are $\overline{\rm O}$-compatible pair
(resp. $\underline{\rm O}$-compatible pair) in both ordered metric
spaces $(X^n,\Delta_n,\sqsubseteq_{n})$ and
$(X^n,\nabla_n,\sqsubseteq_{n})$,
\item [{(iii)}] If $g$ is $\overline{\rm O}$-continuous (resp. $\underline{\rm O}$-continuous) then $G$ is $\overline{\rm O}$-continuous
(resp. $\underline{\rm O}$-continuous) in both ordered metric spaces $(X^n,\Delta_n,\sqsubseteq_{n})$
and $(X^n,\nabla_n,\sqsubseteq_{n})$,
\item [{(iv)}] If $F$ is $\overline{\rm O}$-continuous (resp. $\underline{\rm O}$-continuous) then
$F_{\ast}$ is $\overline{\rm O}$-continuous (resp. $\underline{\rm
O}$-continuous) in both ordered metric spaces
$(X^n,\Delta_n,\sqsubseteq_{n})$ and
$(X^n,\nabla_n,\sqsubseteq_{n})$,
\item [{(v)}] If $F$ is $(g,\overline{\rm O})$-continuous (resp. $(g,\underline{\rm O})$-continuous) then
$F_{\ast}$ is $(G,\overline{\rm O})$-continuous (resp.
$(G,\underline{\rm O})$-continuous) in both ordered metric spaces
$(X^n,\Delta_n,\sqsubseteq_{n})$ and
$(X^n,\nabla_n,\sqsubseteq_{n})$,
\item [{(vi)}] If $(E,d,\preceq)$ has {\it g-ICU} property (resp. {\it g-DCL} property) then both $(E^n,\Delta_n,\sqsubseteq_{n})$ and
$(E^n,\nabla_n,\sqsubseteq_{n})$ have {\it G-ICU} property (resp.
{\it G-DCL} property),
\item [{(vii)}] If $(E,d,\preceq)$ has {\it ICU} property (resp. {\it DCL} property) then both $(E^n,\Delta_n,\sqsubseteq_{n})$ and
$(E^n,\nabla_n,\sqsubseteq_{n})$ have {\it ICU} property (resp. {\it
DCL} property).
\end{enumerate}
\noindent{\bf Proof.} We prove above conclusions only for
$\overline{\rm O}$-analogous and only for the ordered metric space
$(E^n,\Delta_n,\sqsubseteq_{n})$. Their $\underline{\rm
O}$-analogous can analogously be proved. In the similar manner, one
can prove same arguments in the framework of ordered metric space
$(X^n,\nabla_n,\sqsubseteq_{n})$.\\

(i) Let $\{{\rm U}^{(m)}\}$ be an
increasing Cauchy sequence in $(E^n,\Delta_n,\sqsubseteq_{n})$.
Denote U$^{(m)}=(x^{(m)}_1,x^{(m)}_2,...,x^{(m)}_n)$, then by
Proposition 4, each of
$\{x_1^{(m)}\}$,$\{x_2^{(m)}\}$,...,$\{x_n^{(m)}\}$ is an increasing
Cauchy sequence in $(E,d,\preceq)$. By $\overline{\rm
O}$-completeness of $(E,d,\preceq)$, $\exists~x_1,x_2,...,x_n\in E$
such that
$$x_i^{(m)}\stackrel{d}{\longrightarrow} x_i\;{\rm for~each}~i\in
I_n,$$ which using Proposition 3, implies that $${\rm
U}^{(m)}\stackrel{\Delta_n}{\longrightarrow} {\rm U},$$ where ${\rm
U}=(x_1,x_2,...,x_n)$. It follows that
$(E^n,\Delta_n,\sqsubseteq_{n})$ is $\overline{\rm O}$-complete.

(ii) Take a sequence $\{{\rm U}^{(m)}\}\subset X^n$ such that
$\{G{\rm U}^{(m)}\}$ and $\{F_\ast{\rm U}^{(m)}\}$ are increasing
(w.r.t. partial order $\sqsubseteq_{n}$) and
$$G({\rm U}^{(m)})\stackrel{\Delta_n}{\longrightarrow}{\rm W}\;{\rm and}\;F_\ast({\rm U}^{(m)})\stackrel{\Delta_n}{\longrightarrow}{\rm W},$$ for some W$\in X^n$. Write
U$^{(m)}=(x^{(m)}_1,x^{(m)}_2,...,x^{(m)}_n)$ and
W$=(z_1,z_2,...,z_n)$. Then, by using Propositions 3 and 4, we
obtain
$$g(x_i^{(m)})\uparrow z_i\;{\rm and}\;F(x_{i_1}^{(m)},x_{i_2}^{(m)},...,x_{i_n}^{(m)})\uparrow z_i\;\;{\rm for~each}~i\in I_n.\eqno(3)$$
On using (3) and $(\ast,\overline{\rm O})$-compatibility of the pair
$(F,g)$, we have
$$\lim\limits_{m\to \infty}d(gF(x_{i_1}^{(m)},x_{i_2}^{(m)},...,x_{i_n}^{(m)}),F(gx_{i_1}^{(m)},gx_{i_2}^{(m)},...,gx_{i_n}^{(m)}))
=0\;\;{\rm for~each}~i\in I_n$$
$i.e.$
$$\lim\limits_{m\to \infty}d(g(F{\rm U}_i^{(m)\ast}),F(G{\rm U}_i^{(m)\ast}))=0\;\;{\rm for~each}~i\in I_n.\eqno(4)$$
Now, owing to (4), we have
\begin{eqnarray*}
\Delta_n(GF_\ast{\rm U}^{(m)},F_\ast G{\rm U}^{(m)})&=&\frac{1}{n}\sum\limits_{i=1}^{n}d(g(F{\rm U}_i^{(m)\ast}),F(G{\rm U}_i^{(m)\ast}))\\
&\rightarrow& 0\;{\rm as}\;n\rightarrow \infty.
\end{eqnarray*}
It follows that $(F_{\ast},G)$ is  $\overline{\rm O}$-compatible
pair in ordered metric space $(X^n,\Delta_n,\sqsubseteq_{n})$.

The procedure of the proofs of parts (iii) and (iv) are similar to
Lemma 5 and the part (v) and hence is left for readers as an
exercise.\\

(v) Take a sequence $\{{\rm U}^{(m)}\}\subset X^n$ and a ${\rm U}\in
X^n$ such that $\{G{\rm U}^{(m)}\}$ is increasing (w.r.t. partial
order $\sqsubseteq_{n}$) and
$$G({\rm U}^{(m)})\stackrel{\Delta_n}{\longrightarrow} G({\rm U}).$$
Write U$^{(m)}=(x^{(m)}_1,x^{(m)}_2,...,x^{(m)}_n)$ and
U$=(x_1,x_2,...,x_n)$. Then, by using Propositions 3 and 4, we
obtain
$$g(x_i^{(m)})\uparrow  g(x_i)\;{\rm for~each}~i\in
I_n.$$ It follows for each $i\in I_n$ that
$$g(x_{i_1}^{(m)})\uparrow g(x_{i_1}),g(x_{i_2}^{(m)})\uparrow g(x_{i_2}),..., g(x_{i_n}^{(m)})\uparrow g(x_{i_n}). \eqno(5)$$ Using (5) and
$(g,\overline{\rm O})$-continuity of $F$, we get
$$F(x_{i_1}^{(m)},x_{i_2}^{(m)},...,x_{i_n}^{(m)})\stackrel{d}{\longrightarrow}F(x_{i_1},x_{i_2},...,x_{i_n})$$
so that $$F({\rm U}^{(m)\ast}_i)\stackrel{d}{\longrightarrow}F(\rm
U)\;{\rm for~each}~i\in I_n,$$ which, by using Proposition 3 gives
rise
$$F_\ast({\rm U}^{(m)})\stackrel{\Delta_n}{\longrightarrow}F_\ast({\rm U}).$$
Hence, $F_{\ast}$ is $(G,\overline{\rm O})$-continuous in ordered
metric space $(X^n,\Delta_n,\sqsubseteq_{n})$.\\

(vi) Suppose that $(E,d,\preceq)$ has {\it g-ICU} property. Take a
sequence $\{{\rm U}^{(m)}\}\subset E^n$ and a ${\rm U}\in E^n$ such
that $\{{\rm U}^{(m)}\}$ is increasing (w.r.t. partial order
$\sqsubseteq_{n}$) and
$${\rm U}^{(m)}\stackrel{\Delta_n}{\longrightarrow} {\rm U}.$$
Write U$^{(m)}=(x^{(m)}_1,x^{(m)}_2,...,x^{(m)}_n)$ and
U$=(x_1,x_2,...,x_n)$. Then, by Propositions 3 and 4, we obtain
$$x^{(m)}_i\uparrow x_i\;{\rm {for\; each}
}\; i\in I_n,$$ which on using {\it g-ICU} property of
$(E,d,\preceq)$, gives rise
$$g(x^{(m)}_i)\preceq g(x_i)\;{\rm {for\; each}
}\; i\in I_n,$$ or equivalently,
$${\rm U}^{(m)}\sqsubseteq_{\iota_n}{\rm U}.$$
It follows that $(E^n,\Delta_n,\sqsubseteq_{n})$ has {\it G-ICU} property.\\

Analogously, it can be proved that if $(E,d,\preceq)$ has {\it g-DCL} property, then $(E^n,\Delta_n,\sqsubseteq_{n})$ has {\it G-DCL} property.\\

(vii) This result is directly follows from (vi) by setting $g=I,$
the identity mapping.\\

\section{Multi-tupled Coincidence Theorems for Compatible Mappings}
\label{SC:Multi-tupled Coincidence Theorems for Compatible Mappings}
In this section, we prove the results regarding the existence and
uniqueness of $\ast$-coincidence points in ordered metric spaces for compatible pair of mappings.\\

\noindent{\bf Theorem 1.} Let $(X,d,\preceq)$ be an ordered metric
space, $E$ an $\overline{\rm O}$-complete subspace of $X$ and
$\ast\in\mathfrak{B_{n}}$. Let $F:X^{n}\rightarrow X$ and
$g:X\to X $ be two mappings. Suppose that the following conditions
hold:
\begin{enumerate}
\item [{(i)}] $F(X^n)\subseteq g(X)\cap E$,
\item [{(ii)}] $F$ has $g$-monotone property,
\item [{(iii)}] $F$ and $g$ are $(\ast,\overline{\rm
O})$-compatible,
\item [{(iv)}] $g$ is $\overline{\rm O}$-continuous,
\item [{(v)}] either $F$ is $\overline{\rm O}$-continuous or $(E,d,\preceq)$ has {\it g-ICU}
property,
\item [{(vi)}] there exist $x^{(0)}_1,x^{(0)}_2,...,x^{(0)}_n
\in X$ such that$$ g(x^{(0)}_{i}) \preceq
F(x^{(0)}_{i_1},x^{(0)}_{i_2},...,x^{(0)}_{i_n})\;{\rm for~ each}~
i\in I_n,$$
\item [{(vii)}] there
exists $\varphi\in \Omega$ such that
$$\frac{1}{n}\sum\limits_{i=1}^{n}d(F(x_{i_1},x_{i_2},...,x_{i_n}),F(y_{i_1},y_{i_2},...,y_{i_n}))\leq\varphi\Big(\frac{1}{n}\sum\limits_{i=1}^{n}d(gx_i,gy_i)\Big)$$
for all $x_{1},x_{2},...,x_{n},y_{1},y_{2},...,y_{n}\in X$ with\\
\indent\hspace{5mm}$g(x_i)\preceq g(y_i)$ for each $i\in I_n$ or
$g(x_i)\succeq g(y_i)$ for each $i\in I_n$,
\end{enumerate}
or alternately
\begin{enumerate}
\item [{(vii$^\prime$)}] there
exists $\varphi\in \Omega$ such that
$$\max\limits_{i\in I_n}d(F(x_{i_1},x_{i_2},...,x_{i_n}),F(y_{i_1},y_{i_2},...,y_{i_n}))\leq\varphi\Big(\max\limits_{i\in I_n}d(gx_i,gy_i)\Big)$$
for all $x_{1},x_{2},...,x_{n},y_{1},y_{2},...,y_{n}\in X$ with\\
\indent\hspace{5mm}$g(x_i)\preceq g(y_i)$ for each $i\in I_n$ or
$g(x_i)\succeq g(y_i)$ for each $i\in I_n$.
\end{enumerate}
Then $F$ and $g$ have an $\ast$-coincidence point.\\

\noindent{\bf Proof.} We can induce two metrics $\Delta_n$ and
$\nabla_n$, patrial order $\sqsubseteq_n$ and two
self-mappings $F_\ast$ and $G$ on $X^n$ defined as in section 5.
By item (i) of Lemma 8, both ordered metric subspaces
$(E^n,\Delta_n,\sqsubseteq_{n})$ and
$(E^n,\nabla_n,\sqsubseteq_{n})$ are $\overline{\rm O}$-complete.
Further,
\begin{enumerate}
\item [{(i)}] implies that $F_\ast(X^n)\subseteq G(X^n)\cap E^n$ by item (i) of Lemma 3,
\item [{(ii)}] implies that $F_\ast$ is $G$-increasing in ordered set $(X^n,\sqsubseteq_{n})$ by Lemma 5,
\item [{(iii)}] implies that $F_\ast$ and $G$ are $\overline{\rm O}$-compatible in both $(X^n,\Delta_n,\sqsubseteq_{n})$ and $(X^n,\nabla_n,\sqsubseteq_{n})$ by item (ii) of Lemma 8,
\item [{(iv)}] implies that $G$ is $\overline{\rm O}$-continuous in both $(X^n,\Delta_n,\sqsubseteq_{n})$ and
$(X^n,\nabla_n,\sqsubseteq_{n})$ by item (iii) of Lemma 8,
\item [{(v)}] implies that either $F_\ast$ is $\overline{\rm O}$-continuous in both $(X^n,\Delta_n,\sqsubseteq_{n})$ and $(X^n,\nabla_n,\sqsubseteq_{n})$ or both
$(E^n,\Delta_n,\sqsubseteq_{n})$ and
$(E^n,\nabla_n,\sqsubseteq_{n})$ have $G$-{\it MCB} property by
items (iv) and (vi) of Lemma 8
\item [{(vi)}] is equivalent to $G({\rm
U}^{(0)})\sqsubseteq_{n}F_\ast({\rm U}^{(0)})$ where
U$^{(0)}=(x^{(0)}_1,x^{(0)}_2,...,x^{(0)}_n) \in X^n$,
\item [{(vii)}] means that $\Delta_n(F_\ast{\rm U},F_\ast{\rm V})\leq \varphi(\Delta_n(G{\rm U},G{\rm
V}))$ for all U=$(x_{1},x_{2},...,x_{n})$,
V=$(y_{1},y_{2},...,y_{n})\in X^n$ with U$\sqsubseteq_{n}$V or
U$\sqsupseteq_{n}$V,
\item [{(vii$^\prime$)}] means that $\nabla_n(F_\ast{\rm U},F_\ast{\rm V})\leq \varphi(\nabla_n(G{\rm U},G{\rm
V}))$ for all U=$(x_{1},x_{2},...,x_{n})$,
V=$(y_{1},y_{2},...,y_{n})\in X^n$ with U$\sqsubseteq_{n}$V or
U$\sqsupseteq_{n}$V.
\end{enumerate}
Therefore, the conditions (i)-(vii) of Lemma 1 are satisfied in the
context of ordered metric space $(X^n,\Delta_n,\sqsubseteq_{n})$ or
$(X^n,\nabla_n,\sqsubseteq_{n})$ and two self-mappings $F_\ast$ and
$G$ on $X^n$. Thus, by Lemma 1, $F_\ast$ and $G$ have a coincidence
point, which is a $\ast$-coincidence point of $F$ and
$g$ by item (iii) of Lemma 3.\\

Now, we present a dual result corresponding to Theorem 1.\\

\noindent\textbf{Theorem 2.} Theorem 1 remains true if certain
involved terms namely: $\overline{\rm O}$-complete,
$(\ast,\overline{\rm O})$-compatible, $\overline{\rm
O}$-continuous and {\it g-ICU} property are respectively replaced by
$\underline{\rm O}$-complete, $(\ast,\underline{\rm O})$-compatible, $\underline{\rm O}$-continuous and {\it
g-DCL} property provided the assumption (vi) is replaced by the following (besides retaining the rest of the hypotheses):\\
\indent\hspace{5mm} (vi)$^\prime$ there exists
$x^{(0)}_1,x^{(0)}_2,...,x^{(0)}_n \in X$ such that
$$ g(x^{(0)}_{i}) \succeq
F(x^{(0)}_{i_1},x^{(0)}_{i_2},...,x^{(0)}_{i_n})\;{\rm for~ each}~
i\in I_n.$$

\noindent{\bf Proof.} The procedure of the proof of this result is
analogously followed, point by point, by the lines of the proof of
Theorem 1.\\

Now, combining Theorems 1 and 2 and making use of Remarks 1-6, we
obtain the following result:\\

\noindent\textbf{Theorem 3.} Theorem 1 remains true if certain
involved terms namely: $\overline{\rm O}$-complete,
$(\ast,\overline{\rm O})$-compatible, $\overline{\rm
O}$-continuous and {\it g-ICU} property are respectively replaced by
O-complete, $(\ast,{\rm O})$-compatible, O-continuous and {\it
g-MCB} property provided the assumption (vi) is replaced by the following (besides retaining the rest of the hypotheses):\\
\indent\hspace{5mm} (vi)$^{\prime\prime}$ there exists
$x^{(0)}_1,x^{(0)}_2,...,x^{(0)}_n \in X$ such that
$$ g(x^{(0)}_{i}) \preceq
F(x^{(0)}_{i_1},x^{(0)}_{i_2},...,x^{(0)}_{i_n})\;{\rm for~ each}~
i\in I_n$$ or
$$ g(x^{(0)}_{i}) \succeq
F(x^{(0)}_{i_1},x^{(0)}_{i_2},...,x^{(0)}_{i_n})\;{\rm for~ each}~
i\in I_n.$$

Notice that using Remarks 2, 3, 7 and 8, Theorems 1, 2 and 3 provide
their consequences, in which the $\overline{\rm O}$, \underline{\rm
O} and O analogous of metrical notions can be replaced by their
usual senses.\\

Now, we present some consequences of Theorems 1, 2 and 3.\\

\noindent{\bf Corollary 1.} Theorem 1 (similarly Theorem 2 and
Theorem 3) remains true if we replace the condition (vii) by the
following condition:
\begin{enumerate}
\item [{(vii)}$^\prime$] there exists $\varphi\in \Omega$ such that
$$d(F(x_1,x_2,...,x_n),F(y_1,y_2,...,y_n))\leq\varphi\Big(\frac{1}{n}\sum\limits_{i=1}^{n}d(gx_i,gy_i)\Big)$$
for all $x_{1},x_{2},...,x_{n},y_{1},y_{2},...,y_{n}\in X$ with\\
\indent\hspace{5mm}$g(x_i)\preceq g(y_i)$ for each $i\in I_n$ or
$g(x_i)\succeq g(y_i)$ for each $i\in I_n$\\
provided that $\ast$ is permuted.
\end{enumerate}

\noindent{\bf Proof.} Set U=$(x_{1},x_{2},...,x_{n})$,
V=$(y_{1},y_{2},...,y_{n})$ then we have $G({\rm
U})\sqsubseteq_{n}G({\rm V})$ or $G({\rm U})\sqsupseteq_{n}G({\rm
V})$. As $G({\rm U})$ and $G({\rm V})$ are comparable, for each
$i\in I_n$, $G({\rm U}^\ast_i)$ and $G({\rm V}^\ast_i)$ are
comparable w.r.t. partial order $\sqsubseteq_{n}$. Applying the
contractivity condition (vii)$^\prime$ on these points and using
Lemma 6, for each $i\in I_n$, we obtain
\begin{eqnarray*}
d(F(x_{i_1},x_{i_2},...,x_{i_n}),F(y_{i_1},y_{i_2},...,y_{i_n}))&\leq& \varphi\Big(\frac{1}{n}\sum\limits_{k=1}^{n}d(gx_{i_k},gy_{i_k})\Big)\\
&=&
\varphi\Big(\frac{1}{n}\sum\limits_{j=1}^{n}d(gx_j,gy_j)\Big)\;{\rm
as\;} \ast\;{\rm is\; permuted}
\end{eqnarray*}
so that
$$d(F(x_{i_1},x_{i_2},...,x_{i_n}),F(y_{i_1},y_{i_2},...,y_{i_n}))\leq\varphi\Big(\frac{1}{n}\sum\limits_{j=1}^{n}d(gx_j,gy_j)\Big)\;{\rm for~each~}i\in I_n.$$
Taking summation over $i\in I_n$ on both the sides of above
inequality, we obtain
$$\sum\limits_{i=1}^{n}d(F(x_{i_1},x_{i_2},...,x_{i_n}),F(y_{i_1},y_{i_2},...,y_{i_n}))\leq n\varphi\Big(\frac{1}{n}\sum\limits_{j=1}^{n}d(gx_j,gy_j)\Big)$$
so that
$$\frac{1}{n}\sum\limits_{i=1}^{n}d(F(x_{i_1},x_{i_2},...,x_{i_n}),F(y_{i_1},y_{i_2},...,y_{i_n}))\leq\varphi\Big(\frac{1}{n}\sum\limits_{j=1}^{n}d(gx_j,gy_j)\Big)$$
for all $x_{1},x_{2},...,x_{n},y_{1},y_{2},...,y_{n}\in X$ with
$g(x_i)\preceq g(y_i)$ for each $i\in I_n$ or
$g(x_i)\succeq g(y_i)$ for each $i\in I_n$.\\

Therefore, the contractivity condition (vii) of Theorem 1
(similarly Theorem 2 and Theorem 3) holds and hence Theorem 1
(similarly Theorem 2 and Theorem 3) is
applicable.\\

\noindent{\bf Corollary 2.} Theorem 1 (similarly Theorem 2 and
Theorem 3) remains true if we replace the condition (vii$^\prime$)
by the following condition:
\begin{enumerate}
\item [{(vii$^\prime$)}$^\prime$] there exists $\varphi\in \Omega$ such that
$$d(F(x_1,x_2,...,x_n),F(y_1,y_2,...,y_n))\leq\varphi\Big(\max\limits_{i\in I_n}d(gx_i,gy_i)\Big)$$
for all $x_{1},x_{2},...,x_{n},y_{1},y_{2},...,y_{n}\in X$ with\\
\indent\hspace{5mm}$g(x_i)\preceq g(y_i)$ for each $i\in I_n$ or
$g(x_i)\succeq g(y_i)$ for each $i\in I_n$\\
provided that either $\ast$ is permuted or $\varphi$ is increasing
on $[0,\infty)$.
\end{enumerate}

\noindent{\bf Proof.} Set U=$(x_{1},x_{2},...,x_{n})$,
V=$(y_{1},y_{2},...,y_{n})$ then similar to previous corollary, for
each $i\in I_n$, $G({\rm U}^\ast_i)$ and $G({\rm V}^\ast_i)$ are
comparable w.r.t. partial order $\sqsubseteq_{n}$. Applying the
contractivity condition (vii$^\prime)^{\prime}$ on these points and
using Lemma 6, for each $i\in I_n$, we obtain
$$d(F(x_{i_1},x_{i_2},...,x_{i_n}),F(y_{i_1},y_{i_2},...,y_{i_n}))\leq \varphi\Big(\max\limits_{k\in
I_n}d(gx_{i_k},gy_{i_k})\Big)$$
\indent\hspace{8.8cm}${\begin{cases}=\varphi\Big(\max\limits_{j\in
I_n}d(gx_j,gy_j)\Big)\;{\rm if\;} \ast\;{\rm is\; permuted,}\cr
\hspace{0.0in}\leq\varphi\Big(\max\limits_{j\in
I_n}d(gx_j,gy_j)\Big)\;{\rm if\;}\varphi\;{\rm is\;
inceasing.}\cr\end{cases}}$\\
so that
$$d(F(x_{i_1},x_{i_2},...,x_{i_n}),F(y_{i_1},y_{i_2},...,y_{i_n}))\leq\varphi\Big(\max\limits_{j\in
I_n}d(gx_i,gy_i)\Big)\;{\rm for~each~}i\in I_n.$$ Taking maximum
over $i\in I_n$ on both the sides of above inequality, we obtain
$$\max\limits_{i\in
I_n}d(F(x_{i_1},x_{i_2},...,x_{i_n}),F(y_{i_1},y_{i_2},...,y_{i_n}))\leq\varphi\Big(\max\limits_{j\in
I_n}d(gx_j,gy_j)\Big)$$ for all
$x_{1},x_{2},...,x_{n},y_{1},y_{2},...,y_{n}\in X$ with
$g(x_i)\preceq g(y_i)$ for each $i\in I_n$ or
$g(x_i)\succeq g(y_i)$ for each $i\in I_n$.\\

Therefore, the contractivity condition (vii$^\prime$) of Theorem 1
(similarly Theorem 2 and Theorem 3) holds and hence Theorem 1
(similarly Theorem 2 and Theorem 3) is
applicable.\\

Now, we present multi-tupled coincidence theorems for linear and
generalized
linear contractions.\\

\noindent{\bf Corollary 3.} In addition to the hypotheses (i)-(vi)
of Theorem 1 (similarly Theorem 2 and Theorem 3), suppose that one
of the following conditions holds:
\begin{enumerate}
\item [{(viii)}] there exists $\alpha\in [0,1)$ such that
$$\frac{1}{n}\sum\limits_{i=1}^{n}d(F(x_{i_1},x_{i_2},...,x_{i_n}),F(y_{i_1},y_{i_2},...,y_{i_n}))\leq\frac{\alpha}{n} \sum\limits_{i=1}^{n}d(gx_i,gy_i)$$
for all $x_{1},x_{2},...,x_{n},y_{1},y_{2},...,y_{n}\in X$ with\\
$g(x_i)\preceq g(y_i)$ for each $i\in I_n$ or $g(x_i)\succeq g(y_i)$
for each $i\in I_n$,
\item [{(ix)}] there exists $\alpha\in [0,1)$ such that
$$\max\limits_{i\in I_n}d(F(x_{i_1},x_{i_2},...,x_{i_n}),F(y_{i_1},y_{i_2},...,y_{i_n}))\leq\alpha\max\limits_{i\in I_n}d(gx_i,gy_i)$$
for all $x_{1},x_{2},...,x_{n},y_{1},y_{2},...,y_{n}\in X$ with\\
$g(x_i)\preceq g(y_i)$ for each $i\in I_n$ or $g(x_i)\succeq g(y_i)$
for each $i\in I_n$.
\end{enumerate}
Then $F$ and $g$ have an $\ast$-coincidence point.\\

\noindent{\bf Proof.} On setting $\varphi(t)=\alpha t$ with
$\alpha\in [0,1)$, in Theorem 1 (similarly Theorem 2 and Theorem 3),
we get our
result.\\

\noindent{\bf Corollary 4.} In addition to the hypotheses (i)-(vi)
of Theorem 1 (similarly Theorem 2 and Theorem 3), suppose that one
of the following conditions holds:
\begin{enumerate}
\item [{(x)}] there exists $\alpha\in [0,1)$ such that
$$d(F(x_1,x_2,...,x_n),F(y_1,y_2,...,y_n))\leq\alpha\max\limits_{i\in I_n}d(gx_i,gy_i)$$
for all $x_{1},x_{2},...,x_{n},y_{1},y_{2},...,y_{n}\in X$ with\\
$g(x_i)\preceq g(y_i)$ for each $i\in I_n$ or $g(x_i)\succeq g(y_i)$
for each $i\in I_n$,
\item [{(xi)}] there exists $\alpha_1,\alpha_2,...,\alpha_n\in [0,1)$ with $\sum\limits_{i=1}^n\alpha_i<1$ such that
$$d(F(x_1,x_2,...,x_n),F(y_1,y_2,...,y_n))\leq\sum\limits_{i=1}^{n}\alpha_i d(gx_i,gy_i)$$
for all $x_{1},x_{2},...,x_{n},y_{1},y_{2},...,y_{n}\in X$ with\\
$g(x_i)\preceq g(y_i)$ for each $i\in I_n$ or $g(x_i)\succeq g(y_i)$
for each $i\in I_n$,
\item [{(xii)}] there exists $\alpha\in [0,1)$ such that
$$d(F(x_1,x_2,...,x_n),F(y_1,y_2,...,y_n))\leq\frac{\alpha}{n} \sum\limits_{i=1}^{n}d(gx_i,gy_i)$$
for all $x_{1},x_{2},...,x_{n},y_{1},y_{2},...,y_{n}\in X$ with\\
$g(x_i)\preceq g(y_i)$ for each $i\in I_n$ or $g(x_i)\succeq g(y_i)$
for each $i\in I_n$.
\end{enumerate}
Then $F$ and $g$ have an $\ast$-coincidence point.\\

\noindent{\bf Proof.} Setting $\varphi(t)=\alpha t$ with $\alpha\in
[0,1)$, in Corollary 2, we get the result
corresponding to the contractivity condition (x). Notice that here $\varphi$ is increasing on $[0,\infty)$.\\

To prove the result corresponding to (xi), let
$\beta=\sum\limits_{i=1}^n\alpha_i<1$, then we have
\begin{align*}
d(F(x_1,x_2,...,x_n),F(y_1,y_2,...,y_n))&\leq \sum\limits_{i=1}^n\alpha_id(gx_i,gy_i)\nonumber\\
& \leq \Big(\sum\limits_{i=1}^n\alpha_i\Big) \max\limits_{j\in I_n}d(gx_j,gy_j)\nonumber\\
& =\beta \max\limits_{j\in I_n}d(gx_j,gy_j)\nonumber\\
\end{align*}

so that result follows from the result corresponding to (x).\\

Finally, setting $\alpha_i=\frac{\alpha}{n}$ for all $i\in I_n,$
where $\alpha\in [0,1)$ in (xi), we get the result corresponding to
(xii). Notice that here $\sum\limits_{i=1}^n\alpha_i=\alpha<1$.\\

Now, we present uniqueness result corresponding to Theorem 1 (resp.
Theorem 2 and Theorem 3), which
runs as follows:\\

\noindent{\bf Theorem 4.} In addition to the hypotheses of Theorem 1
(resp. Theorem 2 and Theorem 3), suppose that for every pair
$(x_1,x_2,...,x_n)$, $(y_1,y_2,...,y_n)\in X^n$, there exists
$(z_1,z_2,...,z_n)\in X^n$ such that $(gz_1,gz_2,...,gz_n)$ is
comparable to $(gx_1,gx_2,...,gx_n)$ and $(gy_1,gy_2,...,gy_n)$
w.r.t. partial order $\sqsubseteq_{n}$, then $F$ and $g$ have a unique point
of $\ast$-coincidence, which remains also a unique common
$\ast$-fixed point.\\

\noindent{\bf Proof.} Set U=$(x_{1},x_{2},...,x_{n})$,
V=$(y_{1},y_{2},...,y_{n})$ and W=$(z_{1},z_{2},...,z_{n})$, then by
one of our assumptions $G({\rm W})$ is comparable to $G({\rm U})$
and $G({\rm V})$. Therefore, all the conditions of Lemma 1 are
satisfied. Hence, by Lemma 1, $F_\ast$ and $G$ have a unique common
fixed point, a unique point
of coincidence as well as a unique common fixed point, which is
indeed a unique point of $\ast$-coincidence as well as a unique
common $\ast$-fixed point of
$F$ and $g$ by items (iv) and (v) of Lemma 3.\\

\noindent{\bf Theorem 5.} In addition to the hypotheses of Theorem
4, suppose that $g$ is one-one, then $F$ and $g$ have
a unique $\ast$-coincidence point.\\

\noindent{\bf Proof.} Let U=$(x_1,x_2,...,x_n)$ and
V=$(y_1,y_2,...,y_n)$ be two $\ast$-coincidence point of $F$ and $g$
then then using Theorem 4, we obtain
$$(gx_1,gx_2,...,gx_n)=(gy_1,gy_2,...,gy_n)$$
or equivalently
$$g(x_i)=g(y_i)\;{\rm for~each~}i\in I_n.$$
As $g$ is one-one, we have
$$x_i=y_i\;{\rm for~each~}i\in I_n.$$
It follows that U=V, $i.e.$, $F$ and $g$ have a unique $\ast$-coincidence
point.\\

\section{Multi-tupled Coincidence Theorems without Compatibility of mappings}
\label{SC:Multi-tupled Coincidence Theorems without Compatibility of
mappings}

In this section, we prove the results regarding the existence and
uniqueness of $\ast$-coincidence points in an ordered metric space
$X$ for a pair of mappings $F:X^{n}\rightarrow X$ and $g:X\to X$,
which are not necessarily compatible.\\

\noindent{\bf Theorem 6.} Let $(X,d,\preceq)$ be an ordered metric
space, $E$ an $\overline{\rm O}$-complete subspace of $X$ and
$\ast\in\mathfrak{B_{n}}$. Let $F:X^{n}\rightarrow X$ and
$g:X\to X $ be two mappings. Suppose that the following conditions
hold:
\begin{enumerate}
\item [{(i)}] $F(X^n)\subseteq E\subseteq g(X)$,
\item [{(ii)}] $F$ has $g$-monotone property,
\item [{(iii)}] either $F$ is $(g,\overline{\rm O})$-continuous or $F$ and $g$ are continuous or $(E,d,\preceq)$ has {\it ICU} property,
\item [{(iv)}] there exist $x^{(0)}_1,x^{(0)}_2,...,x^{(0)}_n
\in X$ such that$$ g(x^{(0)}_{i}) \preceq
F(x^{(0)}_{i_1},x^{(0)}_{i_2},...,x^{(0)}_{i_n})\;{\rm for~ each}~
i\in I_n,$$
\item [{(v)}] there
exists $\varphi\in \Omega$ such that
$$\frac{1}{n}\sum\limits_{i=1}^{n}d(F(x_{i_1},x_{i_2},...,x_{i_n}),F(y_{i_1},y_{i_2},...,y_{i_n}))\leq\varphi\Big(\frac{1}{n}\sum\limits_{i=1}^{n}d(gx_i,gy_i)\Big)$$
for all $x_{1},x_{2},...,x_{n},y_{1},y_{2},...,y_{n}\in X$ with\\
\indent\hspace{5mm}$g(x_i)\preceq g(y_i)$ for each $i\in I_n$ or
$g(x_i)\succeq g(y_i)$ for each $i\in I_n$,
\end{enumerate}
or alternately
\begin{enumerate}
\item [{(v$^\prime$)}] there
exists $\varphi\in \Omega$ such that
$$\max\limits_{i\in I_n}d(F(x_{i_1},x_{i_2},...,x_{i_n}),F(y_{i_1},y_{i_2},...,y_{i_n}))\leq\varphi\Big(\max\limits_{i\in I_n}d(gx_i,gy_i)\Big)$$
for all $x_{1},x_{2},...,x_{n},y_{1},y_{2},...,y_{n}\in X$ with\\
\indent\hspace{5mm}$g(x_i)\preceq g(y_i)$ for each $i\in I_n$ or
$g(x_i)\succeq g(y_i)$ for each $i\in I_n$.
\end{enumerate}
Then $F$ and $g$ have an $\ast$-coincidence point.\\

\noindent{\bf Proof.} We can induce two metrics $\Delta_n$ and
$\nabla_n$, patrial order $\sqsubseteq_{n}$ and two self-mappings
$F_\ast$ and $G$ on $X^n$ defined as in section 5. By
item (i) of Lemma 8, both ordered metric subspaces
$(E^n,\Delta_n,\sqsubseteq_n)$ and
$(E^n,\nabla_n,\sqsubseteq_n)$ are $\overline{\rm O}$-complete. Further,
\begin{enumerate}
\item [{(i)}] implies that $F_\ast(X^n)\subseteq E^n\subseteq G(X^n)$ by item (ii) of Lemma 3,
\item [{(ii)}] implies that $F_\ast$ is $G$-increasing in ordered set $(X^n,\sqsubseteq_{n})$ by Lemma 5,
\item [{(iii)}] implies that either $F_\ast$ is $(G,\overline{\rm O})$-continuous in both $(X^n,\Delta_n,\sqsubseteq_{n})$ and $(X^n,\nabla_n,\sqsubseteq_{n})$ or
$F_\ast$ and $G$ are continuous in both $(X^n,\Delta_n)$ and
$(X^n,\nabla_n)$ or both $(E^n,\Delta_n,\sqsubseteq_{n})$ and
$(E^n,\nabla_n,\sqsubseteq_{n})$ have {\it ICU} property by Lemma 7
and items (v) and (vii) of Lemma 8
\item [{(iv)}] is equivalent to $G({\rm
U}^{(0)})\sqsubseteq_{n}F_\ast({\rm U}^{(0)})$ where
U$^{(0)}=(x^{(0)}_1,x^{(0)}_2,...,x^{(0)}_n) \in X^n$,
\item [{(v)}] means that $\Delta_n(F_\ast{\rm U},F_\ast{\rm V})\leq \varphi(\Delta_n(G{\rm U},G{\rm
V}))$ for all U=$(x_{1},x_{2},...,x_{n})$,
V=$(y_{1},y_{2},...,y_{n})\in X^n$ with $G({\rm U})\sqsubseteq_{n}G({\rm V})$ or $G({\rm U})\sqsupseteq_{n}G({\rm V})$,
\item [{(v$^\prime$)}] means that $\nabla_n(F_\ast{\rm U},F_\ast{\rm V})\leq \varphi(\nabla_n(G{\rm U},G{\rm
V}))$ for all U=$(x_{1},x_{2},...,x_{n})$,
V=$(y_{1},y_{2},...,y_{n})\in X^n$ with $G({\rm U})\sqsubseteq_{n}G({\rm V})$ or $G({\rm U})\sqsupseteq_{n}G({\rm V})$.
\end{enumerate}
Therefore, the conditions (i)-(v) of Lemma 2 are satisfied in the
context of ordered metric space $(X^n,\Delta_n,\sqsubseteq_{n})$ or
$(X^n,\nabla_n,\sqsubseteq_{n})$ and two self-mappings $F_\ast$ and
$G$ on $X^n$. Thus, by Lemma 2, $F_\ast$ and $G$ have a coincidence
point, which is a $\ast$-coincidence point of $F$ and
$g$ by item (ii) of Lemma 3.\\

Now, we present a dual result corresponding to Theorem 6.\\

\noindent\textbf{Theorem 7.} Theorem 6 remains true if certain
involved terms namely: $\overline{\rm O}$-complete,
$(g,\overline{\rm O})$-continuous and {\it ICU} property are
respectively replaced by $\underline{\rm O}$-complete,
$(g,\underline{\rm O})$-continuous and {\it
DCL} property provided the assumption (iv) is replaced by the following (besides retaining the rest of the hypotheses):\\
\indent\hspace{5mm} (iv)$^\prime$ there exists
$x^{(0)}_1,x^{(0)}_2,...,x^{(0)}_n \in X$ such that
$$ g(x^{(0)}_{i}) \succeq
F(x^{(0)}_{i_1},x^{(0)}_{i_2},...,x^{(0)}_{i_n})\;{\rm for~ each}~
i\in I_n.$$

\noindent{\bf Proof.} The procedure of the proof of this result is
analogously followed, point by point, by the lines of the proof of
Theorem 6.\\

Now, combining Theorems 6 and 7 and making use of Remarks 1-6, we
obtain the following result:\\

\noindent\textbf{Theorem 8.} Theorem 6 remains true if certain
involved terms namely: $\overline{\rm O}$-complete,
$(g,\overline{\rm O})$-continuous and {\it ICU} property are
respectively replaced by O-complete, $(g,{\rm O})$-continuous and
{\it MCB} property provided the assumption (iv) is replaced by the following (besides retaining the rest of the hypotheses):\\
\indent\hspace{5mm} (iv)$^{\prime\prime}$ there exists
$x^{(0)}_1,x^{(0)}_2,...,x^{(0)}_n \in X$ such that
$$ g(x^{(0)}_{i}) \preceq
F(x^{(0)}_{i_1},x^{(0)}_{i_2},...,x^{(0)}_{i_n})\;{\rm for~ each}~
i\in I_n$$ or
$$ g(x^{(0)}_{i}) \succeq
F(x^{(0)}_{i_1},x^{(0)}_{i_2},...,x^{(0)}_{i_n})\;{\rm for~ each}~
i\in I_n.$$

Notice that using Remarks 2 and 8, Theorems 6, 7 and 8 provide their
consequences, in which the $\overline{\rm O}$, \underline{\rm O} and
O analogous of metrical notions can be replaced by their
usual senses.\\

Similar to Corollaries 1-4, the following consequences of
Theorems 5, 6 and 7 hold.\\

\noindent{\bf Corollary 5.} Theorem 6 (similarly Theorem 7 or
Theorem 8) remains true if we replace the condition (v) by the
following condition:
\begin{enumerate}
\item [{(v)}$^\prime$] there exists $\varphi\in \Omega$ such that
$$d(F(x_1,x_2,...,x_n),F(y_1,y_2,...,y_n))\leq\varphi\Big(\frac{1}{n}\sum\limits_{i=1}^{n}d(gx_i,gy_i)\Big)$$
for all $x_{1},x_{2},...,x_{n},y_{1},y_{2},...,y_{n}\in X$ with\\
\indent\hspace{5mm}$g(x_i)\preceq g(y_i)$ for each $i\in I_n$ or
$g(x_i)\succeq g(y_i)$ for each $i\in I_n$\\
provided that $\ast$ is permuted.\\
\end{enumerate}

\noindent{\bf Corollary 6.} Theorem 6 (similarly Theorem 7 or
Theorem 8) remains true if we replace the condition (v$^\prime$) by
the following condition:
\begin{enumerate}
\item [{(v$^\prime$)}$^\prime$] there exists $\varphi\in \Omega$ such that
$$d(F(x_1,x_2,...,x_n),F(y_1,y_2,...,y_n))\leq\varphi\Big(\max\limits_{i\in I_n}d(gx_i,gy_i)\Big)$$
for all $x_{1},x_{2},...,x_{n},y_{1},y_{2},...,y_{n}\in X$ with\\
\indent\hspace{5mm}$g(x_i)\preceq g(y_i)$ for each $i\in I_n$ or
$g(x_i)\succeq g(y_i)$ for each $i\in I_n$\\
provided that either $\ast$ is permuted or $\varphi$ is increasing on $[0,\infty)$.\\
\end{enumerate}

\noindent{\bf Corollary 7.} In addition to the hypotheses (i)-(iv) of
Theorem 6 (similarly Theorem 7 or Theorem 8), suppose that one of
the following conditions holds:
\begin{enumerate}
\item [{(vi)}] there exists $\alpha\in [0,1)$ such that
$$\frac{1}{n}\sum\limits_{i=1}^{n}d(F(x_{i_1},x_{i_2},...,x_{i_n}),F(y_{i_1},y_{i_2},...,y_{i_n}))\leq\frac{\alpha}{n} \sum\limits_{i=1}^{n}d(gx_i,gy_i)$$
for all $x_{1},x_{2},...,x_{n},y_{1},y_{2},...,y_{n}\in X$ with\\
\indent\hspace{5mm}$g(x_i)\preceq g(y_i)$ for each $i\in I_n$ or
$g(x_i)\succeq g(y_i)$ for each $i\in I_n$,
\item [{(vii)}] there exists $\alpha\in [0,1)$ such that
$$\max\limits_{i\in I_n}d(F(x_{i_1},x_{i_2},...,x_{i_n}),F(y_{i_1},y_{i_2},...,y_{i_n}))\leq\alpha\max\limits_{i\in I_n}d(gx_i,gy_i)$$
for all $x_{1},x_{2},...,x_{n},y_{1},y_{2},...,y_{n}\in X$ with\\
\indent\hspace{5mm}$g(x_i)\preceq g(y_i)$ for each $i\in I_n$ or
$g(x_i)\succeq g(y_i)$ for each $i\in I_n$.
\end{enumerate}
Then $F$ and $g$ have an $\ast$-coincidence point.\\

\noindent{\bf Corollary 8.} In addition to the hypotheses (i)-(iv) of
Theorem 6 (similarly Theorem 7 or Theorem 8), suppose that one of
the following conditions hold:
\begin{enumerate}
\item [{(viii)}] there exists $\alpha\in [0,1)$ such that
$$d(F(x_1,x_2,...,x_n),F(y_1,y_2,...,y_n))\leq\alpha\max\limits_{i\in I_n}d(gx_i,gy_i)$$
for all $x_{1},x_{2},...,x_{n},y_{1},y_{2},...,y_{n}\in X$ with\\
\indent\hspace{5mm}$g(x_i)\preceq g(y_i)$ for each $i\in I_n$ or
$g(x_i)\succeq g(y_i)$ for each $i\in I_n$,
\item [{(ix)}] there exists $\alpha_1,\alpha_2,...,\alpha_n\in [0,1)$ with $\sum\limits_{i=1}^n\alpha_i<1$ such that
$$d(F(x_1,x_2,...,x_n),F(y_1,y_2,...,y_n))\leq\sum\limits_{i=1}^{n}\alpha_i d(gx_i,gy_i)$$
for all $x_{1},x_{2},...,x_{n},y_{1},y_{2},...,y_{n}\in X$ with\\
\indent\hspace{5mm}$g(x_i)\preceq g(y_i)$ for each $i\in I_n$ or
$g(x_i)\succeq g(y_i)$ for each $i\in I_n$,
\item [{(x)}] there exists $\alpha\in [0,1)$ such that
$$d(F(x_1,x_2,...,x_n),F(y_1,y_2,...,y_n))\leq\frac{\alpha}{n} \sum\limits_{i=1}^{n}d(gx_i,gy_i)$$
for all $x_{1},x_{2},...,x_{n},y_{1},y_{2},...,y_{n}\in X$ with\\
\indent\hspace{5mm}$g(x_i)\preceq g(y_i)$ for each $i\in I_n$ or
$g(x_i)\succeq g(y_i)$ for each $i\in I_n$.
\end{enumerate}
Then $F$ and $g$ have an $\ast$-coincidence point.\\

Now, we present uniqueness results corresponding to Theorems 6, 7
and 8, which
run as follows:\\

\noindent{\bf Theorem 9.} In addition to the hypotheses of Theorem 6
(similarly Theorem 7 or Theorem 8), suppose that for every pair
$(x_1,x_2,...,x_n)$, $(y_1,y_2,...,y_n)\in X^n$, there exists
$(z_1,z_2,...,z_n)\in X^n$ such that $(gz_1,gz_2,...,gz_n)$ is
comparable to $(gx_1,gx_2,...,gx_n)$ and $(gy_1,gy_2,...,gy_n)$
w.r.t. partial order $\sqsubseteq_n$, then $F$ and $g$ have
a unique point of $\ast$-coincidence.\\

\noindent{\bf Proof.} Set U=$(x_{1},x_{2},...,x_{n})$,
V=$(y_{1},y_{2},...,y_{n})$ and W=$(z_{1},z_{2},...,z_{n})$, then by
one of our assumptions $G({\rm W})$ is comparable to $G({\rm U})$
and $G({\rm V})$. Therefore, all the conditions of Lemma 2 are
satisfied. Hence, by Lemma 2, $F_\ast$ and $G$ have a unique point
of coincidence, which is indeed a unique point of $\ast$-coincidence
of $F$ and $g$ by item (iv) of Lemma 3.\\

\noindent{\bf Theorem 10.} In addition to the hypotheses of Theorem
9, suppose that $g$ is one-one, then $F$ and $g$ have
a unique $\ast$-coincidence point.\\

\noindent{\bf Proof.} Let U=$(x_1,x_2,...,x_n)$ and
V=$(y_1,y_2,...,y_n)$ be two $\ast$-coincidence point of $F$ and $g$
then using Theorem 9, we obtain
$$(gx_1,gx_2,...,gx_n)=(gy_1,gy_2,...,gy_n)$$
or equivalently
$$g(x_i)=g(y_i)\;{\rm for~each~}i\in I_n.$$
As $g$ is one-one, we have
$$x_i=y_i\;{\rm for~each~}i\in I_n.$$
It follows U=V, $i.e.$, $F$ and $g$ have a unique $\ast$-coincidence
point.\\

\noindent{\bf Theorem 11.} In addition to the hypotheses of Theorem
9, suppose that $F$ and $g$ are $(\ast,w)$-compatible, then $F$ and
$g$ have a unique common $\ast$-fixed point.\\

{\noindent\bf{Proof.}} Let $(x_{1},x_{2},...,x_{n})$ be a
$\ast$-coincidence point of $F$ and $g$. Write
$F(x_{i_1},x_{i_2},...,x_{i_n})=g(x_i)=\overline{x}_i$ for each
$i\in I_n$. Then, by Proposition 3,
$(\overline{x}_1,\overline{x}_2,...,\overline{x}_n)$ being a point
of $\ast$-coincidence of $F$ and $g$ is also a $\ast$-coincidence
point of $F$ and $g.$ It follows from Theorem 9 that
$$(gx_1,gx_2,...,gx_n)=(g\overline{x}_1,g\overline{x}_2,...,g\overline{x}_n)$$ $i.e.$,
$\overline{x_i}=g(\overline{x_i})$ for each $i\in I_n$, which for
each $i\in I_n$ yields that
$$F(x_{i_1},x_{i_2},...,x_{i_n})=g(\overline{x_i})=\overline{x_i}.$$
Hence, $(\overline{x}_1,\overline{x}_2,...,\overline{x}_n)$ is a
common $\ast$-fixed point of $F$ and $g$. To prove uniqueness,
assume that $(x^*_1,x^*_2,...,x^*_n)$ is another common $\ast$-fixed
point of $F$ and $g$. Then again from Theorem 9,
$$(gx^*_1,gx^*_2,...,gx^*_n)=(g\overline{x}_1,g\overline{x}_2,...,g\overline{x}_n)$$
$i.e.$
$$(x^*_1,x^*_2,...,x^*_n)=(\overline{x}_1,\overline{x}_2,...,\overline{x}_n).$$
This completes the proof.\\

\section{Multi-tupled Fixed Point Theorems}
\label{SC:Multi-tupled Fixed Point Theorems}

On particularizing $g=I$, the identity mapping on $X$, in the
foregoing results contained in Sections 6 and 7, we obtain the
corresponding $\ast$-fixed point results, which run as follows:\\

\noindent{\bf Theorem 12.} Let $(X,d,\preceq)$ be an ordered metric
space, $F:X^{n}\rightarrow X$ a mapping and
$\ast\in\mathfrak{B_{n}}$. Let $E$ be an $\overline{\rm O}$-complete subspace of
$X$ such that $F(X^n)\subseteq E$. Suppose that the following
conditions hold:
\begin{enumerate}
\item [{(i)}]  $F$ has monotone property,
\item [{(ii)}] either $F$ is $\overline{\rm O}$-continuous or $(Y,d,\preceq)$ has {\it ICU} property,
\item [{(iii)}] there exist $x^{(0)}_1,x^{(0)}_2,...,x^{(0)}_n
\in X$ such that $$x^{(0)}_{i} \preceq
F(x^{(0)}_{i_1},x^{(0)}_{i_2},...,x^{(0)}_{i_n})\;{\rm for~ each}~
i\in I_n$$
\item [{(iv)}] there
exists $\varphi\in \Omega$ such that
$$\frac{1}{n}\sum\limits_{i=1}^{n}d(F(x_{i_1},x_{i_2},...,x_{i_n}),F(y_{i_1},y_{i_2},...,y_{i_n}))=\varphi\Big(\frac{1}{n}\sum\limits_{i=1}^{n}d(x_i,y_i)\Big)$$
for all $x_{1},x_{2},...,x_{n},y_{1},y_{2},...,y_{n}\in X$ with\\
\indent\hspace{5mm}$x_i\preceq y_i$ for each $i\in I_n$ or
$x_i\succeq y_i$ for each $i\in I_n$,
\end{enumerate}
or alternately
\begin{enumerate}
\item [{(iv$^\prime$)}] there
exists $\varphi\in \Omega$ such that
$$\max\limits_{i\in I_n}d(F(x_{i_1},x_{i_2},...,x_{i_n}),F(y_{i_1},y_{i_2},...,y_{i_n}))=\varphi\Big(\max\limits_{i\in I_n}d(x_i,y_i)\Big)$$
for all $x_{1},x_{2},...,x_{n},y_{1},y_{2},...,y_{n}\in X$ with\\
\indent\hspace{5mm}$x_i\preceq y_i$ for each $i\in I_n$ or
$x_i\succeq y_i$ for each $i\in I_n$.
\end{enumerate}
Then $F$ has an $\ast$-fixed point.\\

\noindent\textbf{Theorem 13.} Theorem 12 remains true if certain
involved terms namely: $\overline{\rm O}$-complete, $\overline{\rm
O}$-continuous and {\it ICU} property are respectively replaced by
$\underline{\rm O}$-complete, $\underline{\rm O}$-continuous and
{\it DCL} property provided the assumption (iii) is replaced by the following (besides retaining the rest of the hypotheses):\\
\indent\hspace{5mm} (iii)$^\prime$ there exist
$x^{(0)}_1,x^{(0)}_2,...,x^{(0)}_n \in X$ such that
$$x^{(0)}_{i} \succeq
F(x^{(0)}_{i_1},x^{(0)}_{i_2},...,x^{(0)}_{i_n})\;{\rm for~ each}~
i\in I_n.$$

\noindent\textbf{Theorem 14.} Theorem 12 remains true if certain
involved terms namely: $\overline{\rm O}$-complete, $\overline{\rm
O}$-continuous and {\it ICU} property are respectively replaced by
O-complete, O-continuous and
{\it MCB} property provided the assumption (iii) is replaced by the following (besides retaining the rest of the hypotheses):\\
\indent\hspace{5mm} (iii)$^{\prime\prime}$ there exists
$x^{(0)}_1,x^{(0)}_2,...,x^{(0)}_n \in X$ such that
$$ x^{(0)}_{i} \preceq
F(x^{(0)}_{i_1},x^{(0)}_{i_2},...,x^{(0)}_{i_n})\;{\rm for~ each}~
i\in I_n$$ or
$$ x^{(0)}_{i} \succeq
F(x^{(0)}_{i_1},x^{(0)}_{i_2},...,x^{(0)}_{i_n})\;{\rm for~ each}~
i\in I_n.$$

\noindent{\bf Corollary 9.} Theorem 12 (similarly  Theorem 13 or
Theorem 14) remains true if we replace the condition (iv) by the
following condition:
\begin{enumerate}
\item [{(iv)}$^\prime$] there exists $\varphi\in \Omega$ such that
$$d(F(x_1,x_2,...,x_n),F(y_1,y_2,...,y_n))\leq\varphi\Big(\frac{1}{n}\sum\limits_{i=1}^{n}d(x_i,y_i)\Big)$$
for all $x_{1},x_{2},...,x_{n},y_{1},y_{2},...,y_{n}\in X$ with\\
\indent\hspace{5mm}$x_i\preceq y_i$ for each $i\in I_n$ or
$x_i\succeq y_i$ for each $i\in I_n$\\
provided that $\ast$ is permuted.\\
\end{enumerate}

\noindent{\bf Corollary 10.} Theorem 12 (similarly  Theorem 13 or
Theorem 14) remains true if we replace the condition (iv$^\prime$) by
the following condition:
\begin{enumerate}
\item [{(iv$^\prime$)}$^\prime$] there exists $\varphi\in \Omega$ such that
$$d(F(x_1,x_2,...,x_n),F(y_1,y_2,...,y_n))\leq\varphi\Big(\max\limits_{i\in I_n}d(x_i,y_i)\Big)$$
for all $x_{1},x_{2},...,x_{n},y_{1},y_{2},...,y_{n}\in X$ with\\
\indent\hspace{5mm}$x_i\preceq y_i$ for each $i\in I_n$ or
$x_i\succeq y_i$ for each $i\in I_n$\\
provided that either $\ast$ is permuted or $\varphi$ is increasing on $[0,\infty)$.\\
\end{enumerate}

\noindent{\bf Corollary 11.} Theorem 12 (similarly  Theorem 13 or
Theorem 14) remains true if we replace the condition (iv) by the
following condition:
\begin{enumerate}
\item [{(v)}] there exists $\alpha\in [0,1)$ such that
$$\frac{1}{n}\sum\limits_{i=1}^{n}d(F(x_{i_1},x_{i_2},...,x_{i_n}),F(y_{i_1},y_{i_2},...,y_{i_n}))\leq\frac{\alpha}{n} \sum\limits_{i=1}^{n}d(x_i,y_i)$$
for all $x_{1},x_{2},...,x_{n},y_{1},y_{2},...,y_{n}\in X$ with\\
\indent\hspace{5mm}$x_i\preceq y_i$ for each $i\in I_n$ or
$x_i\succeq y_i$ for each $i\in I_n$,
\item [{(vi)}] there exists $\alpha\in [0,1)$ such that
$$\max\limits_{i\in I_n}d(F(x_{i_1},x_{i_2},...,x_{i_n}),F(y_{i_1},y_{i_2},...,y_{i_n}))\leq\alpha\max\limits_{i\in I_n}d(x_i,y_i)$$
for all $x_{1},x_{2},...,x_{n},y_{1},y_{2},...,y_{n}\in X$ with\\
\indent\hspace{5mm}$x_i\preceq y_i$ for each $i\in I_n$ or
$x_i\succeq y_i$ for each $i\in I_n$.\\
\end{enumerate}

\noindent{\bf Corollary 12.} Theorem 12 (similarly  Theorem 13 or
Theorem 14) remains true if we replace the condition (iv) by the
following condition:
\begin{enumerate}
\item [{(vii)}] there exists $\alpha\in [0,1)$ such that
$$d(F(x_1,x_2,...,x_n),F(y_1,y_2,...,y_n))\leq\alpha\max\limits_{i\in I_n}d(x_i,y_i)$$
for all $x_{1},x_{2},...,x_{n},y_{1},y_{2},...,y_{n}\in X$ with\\
\indent\hspace{5mm}$x_i\preceq y_i$ for each $i\in I_n$ or
$x_i\succeq y_i$ for each $i\in I_n$.
\item [{(viii)}] there exist $\alpha_1,\alpha_2,...,\alpha_n\in [0,1)$ with $\sum\limits_{i=1}^n\alpha_i<1$ such that
$$d(F(x_1,x_2,...,x_n),F(y_1,y_2,...,y_n))\leq\sum\limits_{i=1}^{n}\alpha_i d(x_i,y_i)$$
for all $x_{1},x_{2},...,x_{n},y_{1},y_{2},...,y_{n}\in X$ with\\
\indent\hspace{5mm}$x_i\preceq y_i$ for each $i\in I_n$ or
$x_i\succeq y_i$ for each $i\in I_n$.
\item [{(ix)}] there exists $\alpha\in [0,1)$ such that
$$d(F(x_1,x_2,...,x_n),F(y_1,y_2,...,y_n))\leq\frac{\alpha}{n} \sum\limits_{i=1}^{n}d(x_i,y_i)$$
for all $x_{1},x_{2},...,x_{n},y_{1},y_{2},...,y_{n}\in X$ with\\
\indent\hspace{5mm}$x_i\preceq y_i$ for each $i\in I_n$ or
$x_i\succeq y_i$ for each $i\in I_n$.\\
\end{enumerate}

\noindent{\bf Theorem 15.} In addition to the hypotheses of Theorem
12 (similarly  Theorem 13 or Theorem 14), suppose that for every
pair $(x_1,x_2,...,x_n)$, $(y_1,y_2,...,y_n)\in X^n$, there exists
$(z_1,z_2,...,z_n)\in X^n$ such that $(z_1,z_2,...,z_n)$ is
comparable to $(x_1,x_2,...,x_n)$ and $(y_1,y_2,...,y_n)$ w.r.t.
partial order $\sqsubseteq_n$, then $F$ has a unique $\ast$-fixed point.\\

\section{Conclusion}
\label{SC:Conclusion} We have seen that $\ast$-fixed point theorems
proved in Alam $et\;al.$ \cite{unif1} unify all multi-tupled fixed
point theorems involving mixed monotone property. In the similar
manner, the $\ast$-fixed point theorems proved in this paper unify
all multi-tupled fixed point theorems mixed monotone property, which
substantiate the utility of our results. For instance, in the
following lines, we consider
some special cases of our newly proved results:\\

\noindent $\bullet$ On setting $n=2$ and $\ast=\left[\begin{matrix}
1 &2\\
2 &1\\
\end{matrix}\right]$, we obtain sharpened versions of Theorems 3 and 4.\\

\noindent $\bullet$ On setting $n=3$ and $\ast=\left[\begin{matrix}
1 &2 &3\\
2 &1 &3\\
3 &2 &1\\
\end{matrix}\right]$, we obtain sharpened versions of Theorems 1 and 2.\\

\noindent $\bullet$ On setting $n=4$ and $\ast=\left[\begin{matrix}
1 &2 &3 &4\\
1 &4 &3 &2\\
3 &2 &1 &4\\
3 &4 &1 &2\\
\end{matrix}\right]$,
we obtain quartet fixed/coincidence point theorems involving monotone property, which are variants of quartet coincidence point
theorems of Karapinar \cite{Q0}.\\

\noindent $\bullet$ Taking arbitrary $n$ and
$\ast(i,k)=i_k={\begin{cases}i+k-1\;\;\;\;\;\;\;\;\;\;1\leq
k\leq{n-i+1}\cr \hspace{0.0in}i+k-n-1\;\;\;\;{n-i+2}\leq
k\leq{n}\cr\end{cases}}$, we obtain forward cyclic type $n$-tupled
fixed
point results with monotone property.\\

\noindent $\bullet$ Taking arbitrary $n$ and
$\ast(i,k)=i_k={\begin{cases}i-k+1\;\;\;\;\;\;\;\;\;\;1\leq
k\leq{i}\cr \hspace{0.0in}n+i-k+1\;\;\;\;{i+1}\leq
k\leq{n-1}\cr\end{cases}}$, we obtain backward cyclic type
$n$-tupled fixed
point results with monotone property.\\

\noindent $\bullet$ Taking arbitrary $n$ and
$\ast(i,k)=i_k={\begin{cases}i-k+1\;\;\;\;\;\;\;\;\;\;1\leq
k\leq{i}\cr \hspace{0.0in}k-i+1\;\;\;\;{i+1}\leq
k\leq{n}\cr\end{cases}}$, we obtain 1-skew cyclic type $n$-tupled
fixed point results with monotone property.\\

\noindent $\bullet$ Taking arbitrary $n$ and
$\ast(i,k)=i_k={\begin{cases}i+k-1\;\;\;\;\;\;\;\;\;\;1\leq
k\leq{n-i+1}\cr \hspace{0.0in}2n-i-k+1\;\;\;\;{n-i+2}\leq
k\leq{n}\cr\end{cases}}$, we obtain $n$-skew cyclic type $n$-tupled
fixed point results with monotone property.

\end{document}